\documentclass[12pt,smallextended]{svjour3}
\smartqed
\usepackage{lineno}
\modulolinenumbers[5]

\journalname{Journal of XXX}
\usepackage{amssymb,amsmath,amsfonts}
\usepackage{hyperref}
\hypersetup{
    colorlinks=true,
    linkcolor=blue,
    filecolor=blue,
    urlcolor=blue,
    citecolor=cyan,
    pdfstartview=FitW,
    bookmarks=true,
}
\usepackage{array}
\usepackage{bm}
\usepackage{epsfig}
\usepackage{mathptmx}
\usepackage{diagbox}
\usepackage{graphicx}
\usepackage{color}
\usepackage{caption}
\usepackage{wrapfig}
\usepackage{graphics}
\usepackage{enumerate}
\usepackage{amsxtra}
\usepackage{latexsym}
\usepackage{multirow}
\usepackage{slashbox}
\usepackage{subfigure}
\usepackage{epstopdf}
\usepackage{pdfpages}
\usepackage{algorithm}
\usepackage{algorithmic}

\usepackage{url}
\newcommand{\thmlist}{
\begin{list}{Step 1}
{\setlength{\leftmargin}{0.6 in}\setlength{\labelwidth} {0.5 in}}}

\newcommand{\alglist}{
\begin{list}{Step 1}
{\setlength{\leftmargin}{1.1 in} \setlength{\labelwidth}{1.0 in}}}

 \renewcommand{\proof} {\noindent {\bf Proof.} \quad}
 \newcommand{\eproof} {$\quad \square$}
 
 \newtheorem{assumption}{Assumption}

\renewcommand{\subtitle}[1]{\color{blue}}

\def\blue#1{\color{blue}{#1}\color{black}}

\begin{document}


\title{The regularization continuation method with an adaptive time step control
for linearly constrained optimization problems}
\titlerunning{The regularization continuation method}
\author{Xin-long Luo\textsuperscript{$\ast$}  \and Hang Xiao}
\authorrunning{Luo \and Xiao}

\institute{
     Xin-long Luo
     \at
     Corresponding author. School of Artificial Intelligence, \\
     Beijing University of Posts and Telecommunications, P. O. Box 101, \\
     Xitucheng Road  No. 10, Haidian District, 100876, Beijing China\\
     \email{luoxinlong@bupt.edu.cn}            
     \and
     Hang Xiao
     \at
     School of Artificial Intelligence, \\
     Beijing University of Posts and Telecommunications, P. O. Box 101, \\
     Xitucheng Road  No. 10, Haidian District, 100876, Beijing China \\
     \email{xiaohang0210@bupt.edu.cn}
}

\date{Received: date / Accepted: date}
\maketitle

\begin{abstract}
This paper considers the regularization continuation method and the trust-region
updating strategy for the optimization problem with linear equality constraints.
The proposed method utilizes the linear conservation law of the regularization
continuation method such that it does not need to compute the correction step 
for preserving the feasibility other than the previous continuation methods and 
the quasi-Newton updating formulas for the linearly constrained optimization 
problem. Moreover, the new method uses the special limited-memory 
Broyden-Fletcher-Goldfarb-Shanno (L-BFGS) formula as the preconditioning technique 
to improve its computational efficiency in the well-posed phase, and it uses the 
inverse of the regularized two-sided projection of the Lagrangian Hessian as 
the pre-conditioner to improve its robustness. Numerical results also show that 
the new method is more robust and faster than the traditional
optimization method such as the alternating direction method of multipliers (ADMM),
the sequential quadratic programming (SQP) method (the built-in subroutine
fmincon.m of the MATLAB2020a environment), and the recent continuation method
(Ptctr). The computational time of the new method is about 1/3 of that of SQP 
(fmincon.m). Finally, the global convergence analysis of the new method is also given.
\end{abstract}



\vskip 2mm




\section{Introduction} \label{SUBINT}


In this article, we consider the optimization problem with linear equality
constraints as follows:
\begin{align}
  &\min_{x \in \Re^n} \; f(x)  \nonumber \\
  &\text{subject to} \; \; Ax = b,   \label{LEQOPT}
\end{align}
where $A \in \Re^{m \times n}$ is a matrix and $b \in \Re^{m}$ is a vector.
\blue{This problem has many applications in engineering
fields such as the visual-inertial navigation of an unmanned aerial
vehicle maintaining the horizontal flight \cite{CMFO2009,LLS2021},
constrained sparse regression \cite{BF2010}, sparse signal recovery
\cite{FS2009,VLLW2006}, image restoration and de-noising
\cite{FB2010,NWY2010,ST2010}, the Dantzig selector \cite{LPZ2012}, and
support vector machines \cite{FCG2010}. And there are many practical methods
to solve it such as the sequential quadratic programming (SQP) method
\cite{LJYY2019,NW1999}, the penalty function method \cite{FM1990}, feasible
direction methods (see pp. 515-516, \cite{SY2006}),
and the alternating direction method of multipliers (ADMM \cite{BPCE2011})}.

\vskip 2mm

For the constrained optimization problem \eqref{LEQOPT}, the continuation
method \cite{AG2003,CKK2003,Goh2011,KLQCRW2008,Pan1992,Tanabe1980} is another
method other than the traditional optimization method such as SQP, the penalty
function method and ADMM. The advantage of the continuation method over the SQP method
is that the continuation method is capable of finding many local optimal points
of the non-convex optimization problem by following its trajectory, and it is
even possible to find the global optimal solution
\cite{BB1989,Schropp2000,Yamashita1980}. However, the computational efficiency
of the classical continuation method is inferior to that of the traditional
optimization method such as SQP. Recently, the reference \cite{LLS2021} gives
a continuation method with the trusty time-stepping scheme (Ptctr)
for the problem \eqref{LEQOPT} and it is faster than SQP and the penalty method. In order to
improve the computational efficiency and the robustness of the continuation
method for the large-scale optimization problem further, \blue{we consider a special
limited-memory Broyden-Fletcher-Goldfarb-Shanno (L-BFGS) updating formula
\cite{Broyden1970,Fletcher1970,Goldfarb1970,Shanno1970} as the preconditioned
technique in the well-posed phase and use the inverse of the regularized
two-sided projection of the Lagrangian Hessian as the pre-conditioner in the ill-posed phase}.
Moreover, the new method utilizes the linear conservation law of the regularization
method and it does not need to compute the correction step for preserving the feasibility
other than the previous continuation method \cite{LLS2021} and the quasi-Newton
method \cite{NW1999,SY2006}.

\vskip 2mm

The rest of the paper is organized as follows. In section 2, we give the
regularization continuation method with the switching preconditioned technique
and the trust-region updating strategy for the linearly constrained
optimization problem \eqref{LEQOPT}. In section 3, we analyze the global convergence
of this new method. In section 4,
we report some promising numerical results of the new method, in comparison to the
traditional optimization method such as SQP (the built-in subroutine fmincon.m of the
MATLAB2020a environment \cite{MATLAB}), the alternating direction method of multipliers
(ADMM \cite{BPCE2011}, only for convex problems), and the recent continuation method
(Ptctr \cite{LLS2021}) for some large-scale problems. Finally, we give some
discussions and conclusions in section 5.


\section{The adaptive regularization continuation method}


In this section, we give the regularization continuation method with
the switching preconditioned technique and an adaptive time-step control
based on the trust-region updating strategy \cite{CGT2000} for the linearly
constrained optimization problem \eqref{LEQOPT}. Firstly, we consider
the regularized projection Newton flow based on the KKT conditions
of linearly constrained optimization problem. Then, we give the regularization
continuation method with the trust-region updating strategy to follow this
special ordinary differential equations (ODEs). The new method uses a special
L-BFGS updating formula as the preconditioned technique to improve its computational
efficiency in the well-posed phase, and it uses the inverse of the regularized
two-sided projection of the Lagrangian Hessian as the pre-conditioner to
improve its robustness in the ill-posed phase. Finally, we give a preprocessing
method for the infeasible initial point.

\vskip 2mm

\subsection{The regularization projected Newton flow}

\vskip 2mm

For the linearly constrained optimization problem \eqref{LEQOPT},
its optimal solution $x^{\ast}$ needs to satisfy the Karush-Kuhn-Tucker
conditions (p. 328, \cite{NW1999}) as follows:
\begin{align}
  \nabla_{x} L(x, \, \lambda) &= \nabla f(x) + A^{T} \lambda = 0,
    \label{FOKKTG} \\
  Ax - b & = 0,             \label{FOKKTC}
\end{align}
where the Lagrangian function $L(x, \, \lambda)$ is defined by
\begin{align}
  L(x, \, \lambda) = f(x) + \lambda^{T}(Ax-b).
      \label{LAGFUN}
\end{align}
Similarly to the method of the negative gradient flow for the unconstrained
optimization problem \cite{HM1996}, from the first-order necessary conditions
\eqref{FOKKTG}-\eqref{FOKKTC}, we construct a dynamical system of
differential-algebraic equations for problem \eqref{LEQOPT}
\cite{CL2011,LL2010,Luo2012,LLW2013,Schropp2003} as follows:
\begin{align}
    & \frac{dx}{dt} = - \nabla L_{x}(x, \, \lambda)
      = -\left(\nabla f(x) + A^{T} \lambda \right),  \label{DAGF} \\
    & Ax - b = 0.                      \label{LACON}
\end{align}

\vskip 2mm

By differentiating the algebraic constraint \eqref{LACON} with respect to $t$
and substituting it into the differential equation \eqref{DAGF}, we obtain
\begin{align}
  A\frac{dx}{dt} = - A \left(\nabla f(x) + A^{T} \lambda \right)
  = - A \nabla f(x) - AA^{T} \lambda = 0.    \label{DIFALGC}
\end{align}
If we assume that matrix $A$ has full row rank further, from equation
\eqref{DIFALGC}, we obtain
\begin{align}
   \lambda = - \left(AA^{T} \right)^{-1} A \nabla f(x). \label{LAMBDA}
\end{align}
By substituting $\lambda$ of equation \eqref{LAMBDA} into equation \eqref{DAGF},
we obtain the projected gradient flow \cite{Tanabe1980} for the constrained
optimization problem \eqref{LEQOPT} as follows:
\begin{align}
  \frac{dx}{dt} = - \left( I - A^{T} \left(AA^{T}\right)^{-1}A\right)
  \nabla f(x) = - Pg(x), \label{ODGF}
\end{align}
where $g(x) = \nabla f(x)$ and the projection matrix $P$ is defined by
\begin{align}
  P  = I - A^{T} \left(AA^{T}\right)^{-1}A.  \label{PROMAT}
\end{align}

\vskip 2mm

It is not difficult to verify $P^{2} = P$. That is to say, the projection matrix
$P$ is symmetric and its eigenvalues are either 0 or 1. From Theorem 2.3.1 (see p. 73,
\cite{GV2013}), we know that its matrix 2-norm is
\begin{align}
     \|P\| = 1. \label{MATNP}
\end{align}
We denote $P^{+}$ as the Moore-Penrose generalized inverse of the projection
matrix $P$ (see p. 11, \cite{SY2006}). Since the projection matrix $P$ is symmetric
and $P^{2} = P$, it is not difficult to verify
\begin{align}
      P^{+} = P.     \label{GINVP}
\end{align}
Actually, from equation \eqref{GINVP}, we have $P P^{+}P = P(P)P = P = P^{+}$,
$P^{+}P P^{+} = P^{3} = P$, $\left(P^{+}P\right)^{T} = P^{+}P = P$ and
$\left(P P^{+}\right)^{T} = P P^{+} = P$.

\vskip 2mm

Furthermore, from equation \eqref{PROMAT}, we have $AP = 0$. We denote
$\mathcal{N}(A)$ as the null space of $A$. Since the rank of $A$ is $m$,
we know that the rank of $\mathcal{N}(A)$ equals $n-m$ and there
are $n-m$ linearly independent vectors $x_{i} \, (i = 1, \, \ldots, \, n-m)$ to
satisfy $Ax_{i} = 0 \, (i = 1, \, \ldots, \, n-m)$. From equation \eqref{PROMAT},
we know that those $n-m$ linearly independent vectors $x_{i} \, (i = 1, \, \ldots,
\, n-m)$ satisfy $Px_{i} = x_{i} \, (i = 1, \, \ldots, \, n-m)$. That is to say,
the projection matrix $P$ has $n-m$ linearly independent eigenvectors associated
with eigenvalue 1. Consequently, the rank of $P$ is $n-m$. By combining it
with $AP = 0$, we know that $P$ spans the null space of $A$.


\begin{remark}
If $x(t)$ is the solution of the ODE \eqref{ODGF}, it is not difficult to verify
that $x(t)$ satisfies $A (dx/dt) = 0$. That is to say, if the initial point
$x_{0}$ satisfies $Ax_{0} = b$, the solution $x(t)$ of the projected
gradient flow \eqref{ODGF} also satisfies the feasibility $Ax(t) = b, \; \forall t \ge 0$.
This linear conservation property is very useful when we construct a
structure-preserving algorithm \cite{HLW2006,Shampine1998,Shampine1999} to follow
the trajectory of the ODE \eqref{ODGF} to obtain its steady-state solution $x^{\ast}$.
\end{remark}


If we assume that $x(t)$ is the solution of the ODEs \eqref{ODGF}, by using the
property $P^{2} = P$, we obtain
\begin{align}
 \frac{df(x)}{dt} = \left(\nabla f(x)\right)^{T} \frac{dx}{dt}
 = - (\nabla f(x))^{T} P \nabla f(x) =  - g(x)^{T} P^{2} g(x)
 = - \|Pg(x)\|^{2} \le 0.  \nonumber
\end{align}
That is to say, $f(x)$ is monotonically decreasing along the solution curve $x(t)$
of the dynamical system \eqref{ODGF}. Furthermore, the solution $x(t)$ converges
to $x^{\ast}$ when $f(x)$ is lower bounded and $t$ tends to infinity
\cite{HM1996,Schropp2000,Tanabe1980}, where $x^{\ast}$ satisfies the first-order
Karush-Kuhn-Tucker conditions \eqref{FOKKTG}-\eqref{FOKKTC}. Thus, we can follow
the trajectory $x(t)$ of the ODE \eqref{ODGF} to obtain its steady-state solution
$x^{\ast}$, which is also one stationary point of the original optimization problem
\eqref{LEQOPT}.

\vskip 2mm

However, since the Jacobian $P\nabla^{2}f(x)$ of $Pg(x)$ is rank-deficient, we will
confront the numerical difficulties when we use the explicit ODE method to
follow the projected gradient flow \eqref{ODGF} \cite{AP1998,BCP1996,BJ1998}.
In order to mitigate the stiffness of the ODE \eqref{ODGF}, we use the
generalized inverse $(P\nabla^{2} f(x)P)^{+}$ of the two-sided projection
$P\nabla^{2} f(x)P$ of the Lagrangian  Hessian $\nabla^{2}_{xx}L(x, \, \lambda)$
as the pre-conditioner for the ODE \eqref{ODGF}, which
is used similarly to the system of nonlinear equations \cite{LXL2021}, the
unconstrained optimization problem \cite{HM1996,LXLZ2021,LX2022}, the linear programming
problem \cite{LY2021} and the underdetermined system of nonlinear equations \cite{LX2021}.

\vskip 2mm

Firstly, we integrate the ODE \eqref{ODGF} from zero to $t$, then we obtain
\begin{align}
   x(t) = x(t_{0}) - \int_{0}^{t} Pg(x(\tau)) d \tau
    = x(t_{0}) - P \int_{0}^{t} g(x(\tau))d \tau.  \label{INTODGF}
\end{align}
We denote $z(t) = - \int_{0}^{t} g(x(\tau))d \tau$. Thus, from equation
\eqref{INTODGF}, we have
\begin{align}
    x(t) = x(t_{0}) + Pz(t). \label{XTPROZT}
\end{align}
By substituting it into the ODE \eqref{ODGF}, we obtain
\begin{align}
    P\frac{dz(t)}{dt} = - Pg(x(t_{0})+Pz(t)). \label{PRODGF}
\end{align}

\vskip 2mm

Then, by using the generalized inverse $\left(P\nabla^{2}f(x(t_{0})+Pz(t))P\right)^{+}$
of the Jacobian matrix $P\nabla^{2}f(x(t_{0})+Pz(t))P$ as the pre-conditioner for the
ODE \eqref{PRODGF}, we have
\begin{align}
    P\frac{dz(t)}{dt} = - \left(P\nabla^{2}f(x(t_{0})+Pz(t))P\right)^{+}Pg(x(t_{0})+Pz(t)).
    \label{PCODGF}
\end{align}
We reformulate equation \eqref{PCODGF} as
\begin{align}
    \left(P\nabla^{2}f(x(t_{0})+Pz(t))P\right) \frac{dPz(t)}{dt}
    = - Pg(x(t_{0})+Pz(t)),     \label{DAEZGF}
\end{align}
where we use the property $P^{2} = P$. We let $x(t) = Pz(t) + x(t_{0})$ and
substitute it into equation \eqref{DAEZGF}. Then, we obtain the projected Newton
flow for problem \eqref{LEQOPT} as follows:
\begin{align}
    \left(P \nabla^{2}f(x)P\right) \frac{dx(t)}{dt} = - Pg(x). \label{DODGF}
\end{align}

\vskip 2mm

Although the projected Newton flow \eqref{DODGF} mitigates the stiffness of the
ODE such that we can adopt the explicit ODE method to integrate it on the infinite
interval, there are two disadvantages yet. One is that the two-side projection
$P \nabla^{2}f(x)P$ may be not positive semi-definite. Consequently, it can not
ensure that the objective function $f(x)$ is monotonically decreasing along the
solution $x(t)$ of the ODE \eqref{DODGF}. The other is that the solution $x(t)$
of the ODE \eqref{DODGF} may not satisfy the linear conservation law $Adx(t)/dt = 0$.
In order to overcome these two disadvantages, we use the similar regularization
technique of solving the ill-posed problem \cite{Hansen1994,TA1977} for the
projected Newton flow \eqref{DODGF} as follows:
\begin{align}
    \left(\sigma(x) I + P \nabla^{2}f(x)P\right) \frac{dx(t)}{dt} = - Pg(x),
    \label{TRODGF}
\end{align}
where the regularization parameter $\sigma(x)$ satisfies
$\sigma(x) + \mu_{min}\left(P \nabla^{2}f(x)P\right) \ge \sigma_{min} > 0$.
Here, $\mu_{min}(B)$ represents the smallest eigenvalue of matrix $B$.


\begin{remark} \label{RMLCC}
If we assume that $x(t)$ is the solution of the ODE \eqref{TRODGF}, from the
property $AP = 0$, we have
\begin{align}
    A\left(\sigma(x) I + P \nabla^{2}f(x) P\right) \frac{dx(t)}{dt}
    = - APg(x) = 0. \nonumber
\end{align}
Consequently, we obtain $A \sigma(x) dx(t)/dt = 0$. By integrating it, we obtain
$Ax(t) = Ax(t_0) = b$. That is to say, the solution $x(t)$ of the ODE
\eqref{TRODGF} satisfies the linear conservation law $Ax = b$.
\end{remark}


\begin{remark}
From the property $P^{2} = P$ and the ODE \eqref{TRODGF}, we have
\begin{align}
    \left(\sigma(x) P + P\nabla^{2}f(x)P\right)\frac{dx(t)}{dt} = - Pg(x).
    \label{PTRODGF}
\end{align}
By subtracting equation \eqref{PTRODGF} from equation \eqref{TRODGF}, we obtain
\begin{align}
    \sigma(x) P\frac{dx(t)}{dt} - \sigma(x) \frac{dx(t)}{dt} = 0. \nonumber
\end{align}
Namely, when $x(t)$ is the solution of \eqref{TRODGF}, it satisfies
\begin{align}
    P\frac{dx(t)}{dt} = \frac{dx(t)}{dt}. \label{PDXTEDX}
\end{align}
Consequently, from equations \eqref{TRODGF}, \eqref{PDXTEDX} and
$\sigma(x) + \lambda_{min}\left(P \nabla^{2}f(x)P\right) \ge \sigma_{min} > 0$,
we obtain
\begin{align}
     & \frac{df(x(t))}{dt} = (\nabla f(x))^{T}\frac{dx(t)}{dt}
     = (\nabla f(x))^{T}P \frac{dx(t)}{dt} = (Pg(x))^{T}\frac{dx(t)}{dt}
     \nonumber \\
     & \hskip 2mm  =  - (Pg(x))^{T}
     \left(\sigma(x)I + P^{T}\nabla^{2}f(x)P\right)^{-1}(Pg(x)) \le 0. \nonumber
\end{align}
That is to say, $f(x)$ is monotonically decreasing along the solution $x(t)$
of the ODE \eqref{TRODGF}. Furthermore, the solution $x(t)$ converges to $x^{\ast}$
when $f(x)$ is lower bounded and $\|P\nabla^{2}f(x)P\| \le M$
\cite{HM1996,LQQ2004,Schropp2000,Tanabe1980}, where $M$ is a positive constant
and $x^{\ast}$ is the stationary point of the regularized projection
Newton flow \eqref{TRODGF}. Thus, we can follow the trajectory $x(t)$ of the ODE
\eqref{TRODGF} to obtain its stationary point $x^{\ast}$.
\end{remark}


\subsection{The regularization continuation method} \label{SUBTRCM}


The solution curve $x(t)$ of the ODE \eqref{TRODGF} may  not be efficiently
solved by the general ODE method such as backward differentiation formulas
(BDFs, the subroutine ode15s.m of the MATLAB R2020a environment)
\cite{AP1998,BCP1996,BJ1998,JT1995}. Thus, we need to construct
the particular method for this problem. We apply the first-order explicit Euler
method \cite{SGT2003} to the ODE \eqref{TRODGF}, then we obtain the regularized
projection Newton method:
\begin{align}
    & \left(\sigma_{k} I + P \nabla^{2}f(x_{k})P\right) d_{k} = - Pg(x_{k}),
    \label{TRNEWTON} \\
    & \hskip 2mm x_{k+1} = x_{k} + \alpha_{k} d_{k},
    \label{XK1LST}
\end{align}
where $\alpha_k$ is the time step. When $\alpha_{k} = 1$, the regularized
projection Newton method \eqref{TRNEWTON}-\eqref{XK1LST} equals the
Levenberg-Marquardt method \cite{Levenberg1944,LLT2007,LLS2021,Marquardt1963}.

\vskip 2mm

Since the time step $\alpha_{k}$ of the regularized projection Newton method
\eqref{TRNEWTON}-\eqref{XK1LST} is restricted by the numerical stability
\cite{SGT2003}. That is to say, for the linear test equation
$dx/dt = - \lambda x$, its time step $\alpha_{k}$ is restricted by the
stable region $|1-\alpha_{k}{\lambda}/(\sigma_{k} + \lambda)| \le 1$.
Therefore, the large time step can not be adopted in the steady-state phase.
In order to avoid this disadvantage, similarly to the processing technique of
the nonlinear equations \cite{LXL2021,LY2021,LX2021} and the unconstrained
optimization problem \cite{LXLZ2021,LX2022}, we replace $\alpha_{k}$ with  $\Delta t_{k}/(1+\Delta t_{k})$
in equation \eqref{XK1LST} and let $\sigma_{k} = \sigma_{0}/\Delta t_{k}$ in equation
\eqref{TRNEWTON}. Then, we obtain the regularization continuation
method:
\begin{align}
    & B_{k}d_{k} = - Pg(x_{k}), \;  s_{k} = \frac{\Delta t_{k}}{1 + \Delta t_{k}}d_{k},
    \label{TRCM} \\
    &  \hskip 2mm x_{k+1} = x_{k} + s_{k},     \label{XK1}
\end{align}
where $\Delta t_{k}$ is the time step and $B_{k} = \left((\sigma_{0}/\Delta t_{k})I
+ P \nabla^{2}f(x_{k})P\right)$ or its quasi-Newton approximation.

\begin{remark} The time step $\Delta t_k$ of the regularization
continuation method \eqref{TRCM}-\eqref{XK1} is not restricted by the
numerical stability. Therefore, the large time step $\Delta t_{k}$ can be
adopted in the steady-state phase such that the regularization
continuation method \eqref{TRCM}-\eqref{XK1} mimics the projected Newton
method near the stationary point $x^{\ast}$ and it has the fast convergence
rate. The most of all, the new step $\alpha_{k} = \Delta t_{k}/(\Delta t_{k} + 1)$
is favourable to adopt the trust-region updating strategy to adjust the time step
$\Delta t_{k}$ such that the regularization continuation method
\eqref{TRCM}-\eqref{XK1} accurately follows the trajectory of the regularization
flow \eqref{TRODGF} in the transient-state phase and achieves the fast convergence
rate near its stationary point $x^{\ast}$.
\end{remark}

\vskip 2mm

When $B_{k}$ is updated by the BFGS quasi-Newton formula
\cite{Broyden1970,BNY1987,Fletcher1970,Goldfarb1970,Shanno1970} as follows
\begin{align}
    B_{k+1} = B_{k} + \frac{y_{k}y_{k}^{T}}{y_{k}^{T}s_{k}}
    - \frac{B_{k}s_{k}s_{k}^{T}B_{k}}{s_{k}^{T}B_{k}s_{k}}, \; B_{0} = I,
    \label{BFGS}
\end{align}
where $y_{k} = P g(x_{k+1}) - P g(x_{k}), \; s_{k} = x_{k+1} - x_{k}$,
there is an invariance for the transformation matrix $P$ and we state it as
the following lemma \ref{LEMIVBFGS}.

\vskip 2mm

\begin{lemma} \label{LEMIVBFGS}
Assume that $B_{k}$ is updated by the BFGS quasi-Newton formula \eqref{BFGS} and
$s_{k}$ is solved by equation \eqref{TRCM}, then we have $P(B_{k} - I) = B_{k} - I$ and
$Ps_{k} = s_{k}$ for  $k = 0, \, 1, \, 2, \, \ldots $.
\end{lemma}
\proof We prove this property by induction. When $k = 0$, from $P^{2} = P$, we have
$P(B_{0} - I) = 0 = B_{0} - I$ and $Ps_{0} = s_{0}$. We assume that
$P(B_{l} - I) = B_{l} - I$ and $Ps_{l} = s_{l}$ when $k = l$.
Then, when $k = l+1$, from $P^{2} = P$, $Py_{l} = P(Pg(x_{l+1}) - Pg(x_l))
= y_{l}$ and equation \eqref{BFGS}, we have
\begin{align}
    PB_{l+1} & = PB_{l} + \frac{Py_{l} y_{l}^{T}}{y_{l}^{T}s_{l}}
    - \frac{PB_{l}s_{l}s_{l}^{T}B_{l}}{s_{l}^{T}B_{l}s_{l}} \nonumber \\
    & = P + B_{l} - I + \frac{y_{l} y_{l}^{T}}{y_{l}^{T}s_{l}}
     - \frac{(P + B_l - I)s_{l}s_{l}^{T}B_{l}}{s_{l}^{T}B_{l}s_{l}} \nonumber \\
    & = P - I + B_{l} + \frac{y_{l} y_{l}^{T}}{y_{l}^{T}s_{l}}
    - \frac{B_{l}s_{l}^{T} s_{l}B_{l}}{s_{l}^{T}B_{l}s_{l}}
    = P - I + B_{l+1}. \label{BL1IN}
\end{align}
Consequently, we obtain $P (B_{l+1} - I) = B_{l+1} - I$.

\vskip 2mm

From equation \eqref{TRCM} and equation \eqref{BL1IN}, we have
\begin{align}
     & P B_{l+1} d_{l+1} = P(B_{l+1}d_{l+1})
     = - P(Pg(x_{l+1})) = -Pg(x_{l+1}) = B_{l+1}d_{l+1}, \label{PBDL1} \\
     & P B_{l+1} d_{l+1} = (PB_{l+1}) d_{l+1} = (P + B_{l+1} - I)d_{l+1}.
     \label{PDL1}
\end{align}
Consequently, from equations \eqref{PBDL1}-\eqref{PDL1}, we have
$Pd_{l+1} = d_{l+1}$. By combining it with equation \eqref{TRCM}, we obtain
$Ps_{l+1} = s_{l+1}$. Therefore, we know that the conclusion is true by
induction. \eproof

\vskip 2mm

\begin{remark} \label{REMPLS}
From equations \eqref{TRCM}-\eqref{XK1}, Lemma \ref{LEMIVBFGS} and the property
$AP = 0$, it is not difficult to verify $As_{k} = 0$. Thus, if
the initial point $x_{0}$ is feasible, i.e. $Ax_{0} = b$, $x_{k}$ also satisfies
the linear constraint $Ax_{k} = b$. That is to say, the regularization
continuation method \eqref{TRCM}-\eqref{XK1} satisfies the linear conservation
law such that it does not need to compute the correction step for preserving
the linear feasibility other than the previous continuation method and the quasi-Newton
formula \cite{LLS2021} for the linearly constrained optimization problem \eqref{LEQOPT}.
\end{remark}

\subsection{The adaptive step control}


Another issue is how to adaptively adjust the time step $\Delta t_k$
at every iteration. We borrow the adjustment technique of the trust-region radius
from the trust-region method due to its robustness and its fast convergence
rate \cite{CGT2000,Yuan2015}. According to the linear conservation law of the
regularization continuation method \eqref{TRCM}-\eqref{XK1},
$x_{k+1}$ will preserve the feasibility when $Ax_{k} = b$. That is to say,
$x_{k+1}$ satisfies $Ax_{k+1} = b$. Therefore, we use the objective function
$f(x)$ instead of the nonsmooth penalty function $f(x) + \sigma \|Ax-b\|$ as
the merit function. Similarly to the stepping-time scheme of the ODE method for
the unconstrained optimization problem \cite{Higham1999,LLT2007,LLS2021,LXLZ2021},
we also need to construct a local approximation model of $f(x)$ around $x_{k}$.
Here, we adopt the following quadratic function as its approximation model:
\begin{align}
     q_k(s) = f(x_{k}) + s^{T}g_{k} + \frac{1}{2}s^{T}B_{k}s,
     \label{QOAM}
\end{align}
where $g_{k} = \nabla f(x_{k})$ and $B_{k} = (\sigma_{0}/{\Delta t_{k}})I +
P\nabla^{2} f(x_{k})P$ or its quasi-Newton approximation.

\vskip 2mm

In order to save the computational time, from the regularization continuation
method \eqref{TRCM}-\eqref{XK1}, we simplify the quadratic model
$q_k(s_{k}) - q(0)$ as follows:
\begin{align}
    m_{k}(s_{k}) = g_{k}^{T}s_{k}
    - \frac{0.5\Delta t_{k}}{1+\Delta t_{k}}g_{k}^{T}s_{k}
   = \frac{1+0.5\Delta t_{k}}{1+\Delta t_{k}} g_{k}^{T}s_{k}
   \approx q_{k}(s_{k}) - q_{k}(0).  \label{LOAM}
\end{align}
We enlarge or reduce the time step $\Delta t_k$ at every iteration according
to the following ratio:
\begin{align}
    \rho_k = \frac{f(x_k)-f(x_{k}+s_{k})}{m_k(0)-m_k(s_{k})}.
    \label{MRHOK}
\end{align}
A particular adjustment strategy is given as follows:
\begin{align}
     \Delta t_{k+1} =
        \begin{cases}
          \gamma_1 \Delta t_k, &{\text{if} \hskip 1mm 0 \leq \left|1- \rho_k \right| \le \eta_1,}\\
          \Delta t_k, &{\text{else if} \hskip 1mm \eta_1 < \left|1 - \rho_k \right| < \eta_2,}\\
          \gamma_2 \Delta t_k, &{\text{others},}
        \end{cases} \label{ADTK1}
\end{align}
where the constants are selected as
$\eta_1 = 0.25, \; \gamma_1 = 2, \; \eta_2 = 0.75, \; \gamma_2 = 0.5$  according
to our numerical experiments. We accept the trial step $s_{k}$ and let
$x_{k+1} = x_{k}+s_{k}$, when $\rho_{k} \ge \eta_{a}$ and the approximation model
$m_{k}(0) - m_{k}(s_{k})$ satisfies the Armijo sufficient descent condition:
\begin{align}
     m_{k}(0) - m_{k}(s_{k}) \ge \eta_{m} \|s_{k}\| \|p_{g_k}\|,
     \label{AMGESKPGK}
\end{align}
where $\eta_{a}$ and $\eta_{m}$ are the small positive constants such as
$\eta_{a} = \eta_{m} = 1.0\times 10^{-6}$. Otherwise, we discard it and let
$x_{k+1} = x_{k}$.

\vskip 2mm

\begin{remark}
This new time-stepping scheme based on the trust-region updating strategy
has some advantages, in comparison to the traditional line search strategy
\cite{Luo2005}. If we use the line search strategy and the damped projected
Newton method \eqref{TRNEWTON}-\eqref{XK1LST} to solve the projected Newton flow
\eqref{TRODGF}, in order to achieve the fast convergence rate in the
steady-state phase, the time step $\alpha_{k}$ of the damped projected Newton
method is tried from 1 and reduced by half with many times at every iteration.
Since the linear model $f(x_{k}) + g_{k}^{T}s_{k}$ may not approximate $f(x_{k}+s_{k})$
well in the transient-state phase, the time step $\alpha_{k}$ will be small.
Consequently, the line search strategy consumes the unnecessary trial steps in
the transient-state phase. However, the selection scheme of the time step
based on the trust-region strategy \eqref{MRHOK}-\eqref{ADTK1} can overcome
this shortcoming.
\end{remark}

\subsection{The switching preconditioned technique}


For the large-scale problem, the numerical evaluation of the two-sided projection
$P\nabla^{2}f(x_{k})P$ of the Lagrangian Hessian $\nabla^{2}_{xx} L(x, \, \lambda)$
consumes much time. In order to overcome this shortcoming, in the well-posed
phase, we use the limited-memory BFGS quasi-Newton formula (see
\cite{Broyden1970,Fletcher1970,Goldfarb1970,Mascarenhas2004,Shanno1970} or
pp. 222-230, \cite{NW1999}) to approximate the regularized two-sided projection
$\left(\frac{\sigma_{0}}{\Delta t_{k}}I + P\nabla^{2} f(x_{k})P\right)$ of the
regularization continuation method \eqref{TRCM}-\eqref{XK1}.

\vskip 2mm

Recently, Ullah, Sabi'u and Shah \cite{USS2020} give an efficient L-BFGS updating
formula for the system of monotone nonlinear equations. Furthermore, the reference
\cite{LXLZ2021} also tests its efficiency for some unconstrained optimization problems.
Therefore, we adopt the L-BFGS updating formula to approximate
$\left(\frac{\sigma_{0}}{\Delta t_{k}}I + P\nabla^{2} f(x_{k})P\right)$
in the well-posed phase via slightly revising it as
\begin{align}
      B_{k+1} =
         \begin{cases}
             I - \frac{s_{k}s_{k}^{T}}{s_{k}^{T}s_{k}}
             + \frac{y_{k}y_{k}^{T}}{y_{k}^{T}y_{k}}, \;
             \text{if} \; \left|s_{k}^{T}y_{k}\right| > \theta \|s_{k}\|^{2}, \\
             I, \; \text{otherwise},
        \end{cases}
        \label{LBFGSR}
\end{align}
where $s_{k} = x_{k+1} - x_{k}, \; y_{k} = P\nabla f(x_{k+1}) - P\nabla f(x_{k})$
and $\theta$ is a small positive constant such as $\theta = 10^{-6}$.

\vskip 2mm

By using the Sherman-Morrison-Woodburg formula (p. 17, \cite{SY2006}), from
equation \eqref{LBFGSR}, when $\left|y_{k}^{T}s_{k}\right| > \theta \|s_{k}\|^{2}$,
we obtain the inverse of $B_{k+1}$ as follows:
\begin{align}
       B_{k+1}^{-1} = I  - \frac{y_{k}s_{k}^{T}
       + s_{k}y_{k}^{T}}{y_{k}^{T}s_{k}}
       + 2\frac{y_{k}^{T}y_{k}}{(y_{k}^{T}s_{k})^{2}} s_{k}s_{k}^{T}.
      \label{ILBFGS}
\end{align}
The initial matrix $B_{0}$ can be simply selected as an identity matrix.
From equation \eqref{ILBFGS}, it is not difficult to verify
\begin{align}
       B_{k+1}s_{k} = \frac{y_{k}^{T}s_{k}}{y_{k}^{T}y_{k}} y_{k}. \nonumber
\end{align}
That is to say, $B_{k+1}$ satisfies the scaling quasi-Newton property.

\vskip 2mm

The L-BFGS updating formula \eqref{LBFGSR} has some nice properties such as the
symmetric positive definite property and the positive lower bound of its eigenvalues.

\vskip 2mm

\begin{lemma} \label{LEMHLB}
When $\left|s_{k}^{T}y_{k}\right| > \theta \|s_{k}\|^{2}$, $B_{k+1}$ is
symmetric positive definite and its eigenvalues are greater than
$\left(\theta^{2} \|s_{k}\|^{2}\right)/\left(2\|y_{k}\|^{2}\right)$ and less than 2.
Consequently, when $\left|s_{k}^{T}y_{k}\right| > \theta \|s_{k}\|^{2}$, the
eigenvalues of $B_{k+1}^{-1}$ are greater than ${1}/{2}$ and less than
$\frac{2\|y_{k}\|^{2}}{\theta^{2}\|s_{k}\|^{2}}$.
\end{lemma}

\vskip 2mm

\proof (i) For any nonzero vector $z \in \Re^{n}$, from equation \eqref{LBFGSR},
we have
\begin{align}
     & z^{T}B_{k+1}z = \|z\|^{2}
     - \frac{\left(z^{T}s_{k}\right)^{2}}{\|s_{k}\|^{2}}
     + \frac{\left(z^{T}y_{k}\right)^{2}}{\|y_{k}\|^{2}}
     \ge \frac{\left(z^{T}y_{k}\right)^{2}}{\|y_{k}\|^{2}} \ge 0.
     \label{ZBK1ZGE0}
\end{align}
In the first inequality of equation \eqref{ZBK1ZGE0}, we use the Cauchy-Schwartz
inequality $\|z^{T}s_{k}\| \le \|z\|\|s_{k}\|$ and its equality holds if only if
$z = t s_{k}$. Therefore, $B_{k+1}$ is symmetric positive semi-definite.
When $z = t s_{k}$, since $s_{k}^{T}y_{k} \neq 0$, from equation \eqref{ZBK1ZGE0},
we have $z^{T}B_{k+1}z = t^{2}{\left(s_{k}^{T}y_{k}\right)^{2}}/{\|y_{k}\|^{2}}
> 0$. Consequently, $B_{k+1}$ is symmetric positive definite when $s_{k}^{T}y_{k}
\neq 0$.

\vskip 2mm

(ii) It is not difficult to know that there exist at least $(n-2)$ linearly independent
vectors $z_{1}, \, z_{2}, \, \ldots, \, z_{n-2}$ to satisfy $s_{k}^{T}z_{i}= 0, \,
y_{k}^{T}z_{i} = 0 \, ( i = 1 :(n-2))$. That is to say, matrix $B_{k+1}$ defined
by equation \eqref{LBFGSR} has at least $(n-2)$ linearly independent eigenvectors
associated with eigenvalue 1. We denote the other two eigenvalues of
$B_{k+1}$ as $\mu_{i}^{k+1} \, (i = 1:2)$ and set
$\text{tr}(C) = \sum_{i=1}^{n} c_{ii}, \, C \in \Re^{n \times n}$. Then, we have
$\text{tr}(B_{k+1}) = \mu_{1}^{k+1}+\mu_{2}^{k+1}+(n-2)$. By substituting it into
equation \eqref{LBFGSR}, we obtain
\begin{align}
     & \mu_{1}^{k+1} + \mu_{2}^{k+1} = \text{tr}(B_{k+1}) - (n-2) \nonumber \\
     & \hskip 2mm = \text{tr}(I) -
     \text{tr}\left(\frac{s_{k}s_{k}^{T}}{s_{k}^{T}s_{k}}\right)
     + \text{tr}\left(\frac{y_{k}y_{k}^{T}}{y_{k}^{T}y_{k}}\right) - (n-2)
     = 2, \label{TWOEIGSUM}
\end{align}
where we use the property $\text{tr}\left(AB^{T}\right)
= \text{tr}\left(B^{T}A\right)$ of matrices $A, \, B \in \Re^{m \times n}$.
Since matrix $B_{k+1}$ is symmetric positive definite, we know that its
eigenvalues are greater than 0, namely $\mu_{i}^{k+1} > 0 \, (i=1, \, 2)$. By
substituting it into equation \eqref{TWOEIGSUM}, we obtain
\begin{align}
    0 < \mu_{i}^{k+1} < \mu_{1}^{k+1}+\mu_{2}^{k+1} = 2, \; i = 1, \, 2.
    \label{MUILE2}
\end{align}
Furthermore, the symmetric matrix $B_{k+1}$ has a multiple eigenvalue 1 associated
with $(n-2)$ linearly independent eigenvectors. Therefore, by combining it with
equation \eqref{MUILE2}, we know that the eigenvalues of matrix $B_{k+1}$ are
less than 2.

\vskip 2mm

We denote $\mu_{i}^{k+1} \, (i=1:n)$ as the eigenvalues of $B_{k+1}$. Then, we have
$\mu_{i}^{k+1} = 1 \, (i=3:n)$. By using the property
$\det(B_{k+1}) = \prod_{i=1}^{n} \mu_{i}^{k+1} = \mu_{1}^{k+1} \mu_{2}^{k+1}$,
from equation \eqref{LBFGSR}, we obtain
\begin{align}
    & \mu_{1}^{k+1}\mu_{2}^{k+1}
    = \det(B_{k+1}) = \det\left(\left(I + \frac{y_{k}y_{k}^{T}}{y_{k}^{T}y_{k}}\right)
    \left(I - \left(I + \frac{y_{k}y_{k}^{T}}{y_{k}^{T}y_{k}}\right)^{-1}s_{k}
    \frac{s_{k}^{T}}{s_{k}^{T}s_{k}}\right)\right) \nonumber \\
    & \hskip 2mm
    = \det\left(I + \frac{y_{k}y_{k}^{T}}{y_{k}^{T}y_{k}}\right)
    \det\left(I - \left(I + \frac{y_{k}y_{k}^{T}}{y_{k}^{T}y_{k}}\right)^{-1}s_{k}
    \frac{s_{k}^{T}}{s_{k}^{T}s_{k}}\right) \nonumber \\
   & \hskip 2mm
   = 2 \left(1 - \frac{1}{\|s_{k}\|^{2}}s_{k}^{T}
   \left(I + \frac{y_{k}y_{k}^{T}}{y_{k}^{T}y_{k}}\right)^{-1}s_{k} \right)
   \nonumber \\
   & \hskip 2mm
   = 2 \left(1 - \frac{1}{\|s_{k}\|^{2}}s_{k}^{T}
   \left(I- \frac{y_{k}y_{k}^{T}}{2y_{k}^{T}y_{k}}\right)s_{k} \right)
   = \frac{\left(s_{k}^{T}y_{k}\right)^{2}}
   {\left(y_{k}^{T}y_{k}\right)\left(s_{k}^{T}s_{k}\right)}. \label{TWOEIGPROD}
\end{align}
From equation \eqref{MUILE2}, we know $ 0 < \mu_{i}^{k} < 2 \, (i = 1, \, 2)$. By substituting
it into equation \eqref{TWOEIGPROD}, we obtain
\begin{align}
    \mu_{i}^{k+1} > \frac{1}{2}\frac{\left(s_{k}^{T}y_{k}\right)^{2}}
   {\|s_{k}\|^{2} \|y_{k}\|^{2}}, \; i = 1, \, 2.  \label{MU12GE}
\end{align}
By combining it with $\mu_{i}^{k+1} = 1 \, (i = 3:n)$, we have
\begin{align}
    \mu_{i}^{k+1} \ge \min\left\{1, \; \frac{1}{2}\frac{\left(s_{k}^{T}y_{k}\right)^{2}}
   {\|s_{k}\|^{2} \|y_{k}\|^{2}}\right\}
   = \frac{1}{2}\frac{\left(s_{k}^{T}y_{k}\right)^{2}}
   {\|s_{k}\|^{2} \|y_{k}\|^{2}}
   \ge \frac{1}{2}\theta^{2}\frac{\|s_{k}\|^{2}}{\|y_{k}\|^{2}}, \label{MUIGE}
\end{align}
where we use the Cauchy-Schwartz inequality $|s_{k}^{T}y_{k}| \le \|s_{k}\| \|y_{k}\|$.

\vskip 2mm

Since the matrix $B_{k+1}$ is symmetric positive definite when
$\left|s_{k}^{T}y_{k}\right| > \theta \|s_{k}\|^{2}$, the inverse of $B_{k+1}$
exists. Furthermore, the eigenvalues of $B_{k+1}^{-1}$ equal
$1/\mu_{i}^{k+1} \, (i = 1:n)$. By combining it with equations
\eqref{MUILE2} and \eqref{MUIGE}, we know that the eigenvalues of $B_{k+1}^{-1}$
are greater than 1/2 and less than $\left(2\|y_{k}\|^{2}\right)/\left(\theta^{2} \|s_{k}\|^{2}\right)$
when $\left|s_{k}^{T}y_{k}\right| > \theta \|s_{k}\|^{2}$.  \qed

\vskip 2mm

According to our numerical experiments \cite{LXLZ2021}, the L-BFGS
updating formula \eqref{LBFGSR} works well for most problems and the
objective function decreases very fast in the well-posed phase.
However, for the ill-posed problems, the L-BFGS updating formula \eqref{LBFGSR}
will approach the stationary
solution $x^{\ast}$ very slow in the ill-posed phase. Furthermore, it fails
to get close to the stationary solution $x^{\ast}$ sometimes.

\vskip 2mm

In order to improve the robustness of the regularization continuation
method \eqref{TRCM}-\eqref{XK1}, we adopt the inverse $B_{k+1}^{-1}$ of the
regularized two-side projection of the Lagrangian Hessian
$\nabla^{2}_{xx}L(x, \, \lambda)$  as the pre-conditioner in the ill-posed
phase, where $B_{k+1}$ is defined by
\begin{align}
    B_{k+1} = \frac{\sigma_{0}}{\Delta t_{k+1}}I + P\nabla^{2} f(x_{k+1})P.
    \label{RPHM}
\end{align}
Now, the problem is how to automatically identify the ill-posed phase and
switch to the inverse of the regularized two-sided projection from the L-BFGS
updating formula \eqref{LBFGSR}. Here, we adopt the simple switching criterion.
Namely, we regard that the regularization continuation method
\eqref{TRCM}-\eqref{XK1} is in the ill-posed phase once there exists the time
step $\Delta t_{K} \le 10^{-3}$.

\vskip 2mm

In the ill-posed phase,  the computational time of the two-sided projection
of the Lagrangian Hessian $\nabla^{2}_{xx} L(x, \, \lambda)$
is heavy if we update the two-sided projection $P\nabla^{2}f(x_{k})P$
at every iteration. In order to save its computational time, we set
$B_{k+1} = B_{k}$ when $m_{k}(0) - m_{k}(s_{k})$
approximates $f(x_{k}) - f(x_{k}+s_{k})$ well, where the approximation model
$m_{k}(s_{k})$ is defined by equation \eqref{LOAM}. Otherwise, we update
$B_{k+1} = \left(\frac{\sigma_{0}}{\Delta t_{k+1}}I + P\nabla^{2}f(x_{k+1})P\right)$ in the
ill-posed phase. In the ill-posed phase, a practice updating strategy is give by
\begin{align}
     B_{k+1}
       = \begin{cases}
            B_{k},  \; \text{if} \; |1- \rho_{k}| \le \eta_{1}, \\
            \frac{\sigma_{0}}{\Delta t_{k+1}}I + P\nabla^{2}f(x_{k+1})P, \; \text{otherwise},
          \end{cases} \label{UPDJK1}
\end{align}
where $\rho_{k}$ is defined by equations \eqref{LOAM}-\eqref{MRHOK} and
$\eta_{1} = 0.25$.

\vskip 2mm

For a real-world problem, the analytical Hessian matrix  $\nabla^{2}f(x_{k})$ may not be
offered. Thus, in practice, we replace the two-sided projection $P\nabla^{2}f(x_{k})P$
with its difference approximation as follows:
\begin{align}
     P\nabla^{2}f(x_{k})P \approx
    \left[\frac{Pg(x_{k} + \epsilon Pe_{1}) - Pg(x_{k})}{\epsilon}, \,
    \ldots, \, \frac{Pg(x_{k} + \epsilon Pe_{n}) - Pg(x_{k})}{\epsilon}\right],
    \label{NUMHESS}
\end{align}
where the elements of $e_{i}$ equal 0 except for the $i$-th element equaling 1,
and the parameter $\epsilon$ can be selected as $10^{-6}$ according to our numerical
experiments.

\subsection{The treatment of rank-deficient problems and infeasible initial points}
\label{SUBSECDEG}

\vskip 2mm

For a real-world problem, matrix $A$ may be deficient-rank. We assume that the
rank of $A$ is $r$ and we use the QR decomposition (pp.276-278, \cite{GV2013}) to
factor $A^{T}$ into a product of an orthogonal matrix $Q \in \Re^{n \times n}$
and an upper triangular matrix $R \in \Re^{n \times m}$ as follows:
\begin{align}
      A^{T}E = QR = \begin{bmatrix} Q_{1} | Q_{2} \end{bmatrix}
      \begin{bmatrix} R_{1} \\ 0     \end{bmatrix}, \label{ATQR}
\end{align}
where $E \in \Re^{m \times m}$ is a permutation matrix,  $R_{1} = R(1:r, \, 1:m)$
is an upper triangular matrix and its diagonal elements are non-zero, and $Q_{1} = Q(1:n, \, 1:r)$,
$Q_{2} = Q(1:n, \, (r+1):n)$ satisfy $Q_{1}^{T}Q_{1} = I$, $Q_{2}^{T}Q_{2} = I$ and
$Q_{1}^{T}Q_{2} = 0$. Then, we reduce the linear constraint $Ax = b$ to
\begin{align}
     Q_{1}^{T}x = b_{r}, \label{RLEC}
\end{align}
where $b_{r} = \left(R_{1}R_{1}^{T}\right)^{-1}\left(R_{1}\left(E^{T}b\right)\right)$.

\vskip 2mm

From equations \eqref{PROMAT} and \eqref{RLEC}, we simplify the projection matrix
$P$ as
\begin{align}
      P = I - Q_{1}Q_{1}^{T} = Q_{2}Q_{2}^{T}. \label{SMPPROJ}
\end{align}
In practical computation, we adopt the different formulas of the projection
matrix $P$ according to $r \le n/2$ or $ r > n/2$. Thus, we give the computational
formula of the projected gradient $Pg_{k}$ as follows:
\begin{align}
       Pg_{k} = \begin{cases}
                     g_{k} - Q_{1}\left(Q_{1}^{T}g_{k}\right), \;
                     \text{if} \; r \le \frac{1}{2}n, \\
                     Q_{2}\left(Q_{2}^{T}g_{k}\right), \; \text{otherwise}.
                \end{cases}
                \label{PROGK}
\end{align}
where $r$ is the number of columns of $Q_{1}$, i.e. the rank of $A$.

\vskip 2mm

For a real-world optimization problem \eqref{LEQOPT}, we probably meet the
infeasible initial point $x_{0}$. In other words, the initial point may not
satisfy the constraint $Ax = b$. We handle this problem by solving the following
projection problem:
\begin{align}
     \min_{x \in \Re^{n}} \; \left\|x - x_{0} \right\|^{2}
     \; \text{subject to} \hskip 2mm  Q_{1}^{T} x = b_{r},  \label{MINDISTVB}
\end{align}
where $b_{r} = \left(R_{1}R_{1}^{T}\right)^{-1}\left(R_{1}\left(E^{T}b\right)\right)$.
By using the Lagrangian multiplier method to solve problem \eqref{MINDISTVB},
we obtain the initial feasible point $x_{0}^{F}$ of problem \eqref{LEQOPT} as
follows:
\begin{align}
    x_{0}^{F} = x_{0} - Q_{1}\left(Q_{1}^{T}x_{0}-b_{r}\right).
    \label{INFIPT}
\end{align}
For convenience, we set $x_{0} = x_{0}^{F}$ in line 4 of Algorithm \ref{ALGTRCMTR}.

\vskip 2mm

According to the above discussions, we give the detailed implementation of
the regularization continuation method with the trust-region updating
strategy for the linearly constrained optimization problem \eqref{LEQOPT}
in Algorithm \ref{ALGTRCMTR}.

\begin{algorithm}
	\renewcommand{\algorithmicrequire}{\textbf{Input:}}
	\renewcommand{\algorithmicensure}{\textbf{Output:}}
    \newcommand{\algorithmicbreak}{\textbf{break}}
    \newcommand{\BREAK}{\STATE \algorithmicbreak}
	\caption{The regularization continuation method with the trust-region updating strategy
    for linearly constrained optimization problems (Rcmtr)}
    \label{ALGTRCMTR}	
	\begin{algorithmic}[1]
		\REQUIRE ~~ the objective function $f: \; \Re^{n} \to \Re$, the linear constraint
        $Ax  = b, \; A \in \Re^{m \times n}, \; b \in \Re^{m}$,
        the initial point $x_0$ (optional), the tolerance error $\epsilon$ (optional).
		\ENSURE ~~
        the optimal approximation solution $x^{\ast}$.
        \STATE If $x_{0}$ or $\epsilon$ is not provided, we set
        $x_0 = \text{ones}(n, \, 1)$ or $\epsilon = 10^{-6}$.
        \STATE Initialize the parameters: $\eta_{a} = 10^{-6}, \; \eta_{m} = 10^{-10}$,
        $\eta_1 = 0.25, \; \gamma_1 =2, \; \eta_2 = 0.75, \; \gamma_2 = 0.5$,
        $\theta = 10^{-6}$, max\_itc = 300.
        Set $\sigma_{0} = 10^{-4}, \; \Delta t_0 = 10^{-2}$, flag\_illposed\_phase = 0, flag\_success\_trialstep = 1,
        $s_{-1} = 0, \; y_{-1} = 0, \; \rho_{-1} = 0, \; B_{0} = I, \; H_{0} = I$, itc = 0.
        \STATE Factorize matrix $A^{T}$ into $A^{T}E= Q_{1}R_{1}$ with the QR decomposition \eqref{ATQR}.
        Solve the linear system $\left(R_{1}R_{1}^{T}\right)b_{r} = R_{1}\left(E^{T}b\right)$ to obtain $b_{r}$.
        \STATE Compute
        $$ x_{0} \leftarrow x_{0} - Q_{1}\left(Q_{1}^{T}x_{0}-b_{r}\right),
        $$
        such that $x_{0}$ satisfies the linear constraint $Ax_{0} = b$.
        \STATE Set $k = 0$. Evaluate $f_0 = f(x_0)$ and $g_0 = \nabla f(x_0)$.
        \STATE Compute the projected gradient $p_{g_{0}} = Pg_{0}$ according to
        the formula \eqref{PROGK}.
        \WHILE{$\left(\left(\|p_{g_k}\|> \epsilon\right) \, \text{and} \,  (\text{itc} < \text{max\_itc})\right)$}
           \STATE itc = itc + 1;
           \IF{$\Delta t_{k} < 10^{-3}$}
              \STATE Set flag\_illposed\_phase = 1.
           \ENDIF
           \IF{(flag\_illposed\_phase == 0)}
               \IF{(flag\_success\_trialstep == 1)}
                   \IF{$\left(|s_{k-1}^{T}y_{k-1}| > \theta \|s_{k-1}\|^{2}\right)$}
                      \STATE $d_{k} = - \left(p_{g_{k}} - \frac{y_{k-1}(s_{k-1}^{T}p_{g_{k}})
                      + s_{k-1}(y_{k-1}^{T}p_{g_{k}})}{y_{k-1}^{T}s_{k-1}}
                      + 2 \frac{\|y_{k-1}\|^{2}(s_{k-1}^{T}p_{g_{k}})}{(y_{k-1}^{T}s_{k-1})^{2}}s_{k-1}\right)$.
                   \ELSE
                        \STATE $d_{k} = - p_{g_k}$.
                   \ENDIF
               \ENDIF
           \ELSE
               \IF{(flag\_success\_trialstep == 0)}
                  \STATE Set $B_{k} = ({\sigma_{0}}/{\Delta t_{k}})I + H_{k}$ and factorize $B_{k}$ into
                   $B_{k}= Q_{k}R_{k}$ with the QR decomposition.
                \ELSIF{($(|\rho_{k-1} - 1| > 0.25)$}
                   \STATE Evaluate  $H_{k} = P\nabla^{2}f(x_{k})P$ from equation \eqref{NUMHESS}.
                   \STATE Set $B_{k} = ({\sigma_{0}}/{\Delta t_{k}})I + H_{k}$ and factorize $B_{k}$
                   into $B_{k}= Q_{k}R_{k}$ with the QR decomposition.
                \ELSE
                   \STATE $Q_{k} = Q_{k-1}, \; R_{k} = R_{k-1}$.
               \ENDIF
               \STATE Solve the linear system $R_{k}d_{k} = - Q_{k}^{T}p_{g_k}$ to obtain $d_{k}$.
          \ENDIF
          \STATE Set $s_{k} = \frac{\Delta t_{k}}{1+\Delta t_{k}}d_{k}$ and $x_{k+1} = x_{k} + s_{k}$.
          \STATE Evaluate $f_{k+1} = f(x_{k+1})$ and compute the ratio $\rho_{k}$
          from equations \eqref{LOAM}-\eqref{MRHOK}.
          \IF{($\rho_k \ge \eta_{a}$ and $s_{k}$ satisfies the sufficient descent condition \eqref{AMGESKPGK})}
             \STATE Set flag\_success\_trialstep = 1 and evaluate $g_{k+1} = \nabla f(x_{k+1})$.
             \STATE Compute $p_{g_{k+1}} = Pg_{k+1}$ according to the formula
            \eqref{PROGK}. Set $y_{k} = p_{g_{k+1}} - p_{g_{k}}$.
          \ELSE
            \STATE Set flag\_success\_trialstep = 0 and $x_{k+1} = x_{k}, \; f_{k+1} = f_{k}$,
             $p_{g_{k+1}} = p_{g_{k}}, \; g_{k+1} = g_{k}, \; H_{k+1} = H_{k}, \; d_{k+1} = d_{k}.$
          \ENDIF
          \STATE Adjust the time step $\Delta t_{k+1}$ according to the
          trust-region updating strategy \eqref{ADTK1}.
          \STATE Set $k \leftarrow k+1$.
        \ENDWHILE
	\end{algorithmic}
\end{algorithm}

\section{Algorithm Analysis}

In this section, we analyze the global convergence of the regularization
continuation method \eqref{TRCM}-\eqref{XK1} with the trust-region updating
strategy and the switching preconditioned technique for the linearly constrained
optimization problem (i.e. Algorithm \ref{ALGTRCMTR}). Firstly, we give a
lower-bounded estimation of $m_{k}(0) - m_{k}(s_{k})$ $(k = 1, \, 2, \, \ldots)$.
This result is similar to that of the trust-region method for the unconstrained
optimization problem \cite{Powell1975}. For simplicity, we assume that the rank
of matrix $A$ is full and $f$ satisfies Assumption \ref{ASSFUN}.

\vskip 2mm

\begin{assumption} \label{ASSFUN}
Assume that $f(\cdot)$ is twice continuously differential and there exists a
positive constant $M$ such that
\begin{align}
    \left\|\nabla^{2}f(x)\right\| \le M, \label{BOUNDHESS}
\end{align}
holds for all $x \in S_{f}$, where $S_{f} = \{x: \; Ax = b \}$.
\end{assumption}

\vskip 2mm

By combining the property $\|P\| = 1$ of the projection matrix $P$, from the
assumption  \eqref{BOUNDHESS}, we obtain
\begin{align}
    \left\|P\nabla^{2}f(x)P\right\| \le \|P\| \left\|\nabla^{2}f(x)\right\| \|P\|
    = \left\|\nabla^{2}f(x)\right\| \le M. \label{BOUNDPHESS}
\end{align}
According to the property of the matrix norm, we know that the absolute eigenvalue
of $P\nabla^{2} f(x)P$ is less than $M$. We denote $\mu(C)$ as the eigenvalue
of matrix $C$. Then, we know that the eigenvalue of
$((\sigma_{0}/\Delta t)I + P\nabla^{2} f(x) P)$
is $\sigma_{0}/\Delta t + \mu(P\nabla^{2}f(x)P)$. Consequently, from
equation \eqref{BOUNDPHESS}, we known that
\begin{align}
   \frac{\sigma_{0}}{\Delta t}I + P\nabla^{2}f(x)P \succ 0, \;  x \in S_{f}, \;
   \text{when} \; \Delta t < \frac{\sigma_{0}}{M}.    \label{SPDCON}
\end{align}

\vskip 2mm

\begin{lemma} \label{LBSOAM}
Assume that the approximation model $m_{k}(s)$ is defined by equation \eqref{LOAM}
and $s_{k}$ is computed by the regularization continuation method
\eqref{TRCM}-\eqref{XK1}, where matrices $B_{k} \, (k = 1, \, 2, \, \ldots)$
are updated by the L-BFGS formula \eqref{LBFGSR} in the well-posed phase.
Then, we have
\begin{align}
    m_{k}(0) - m_{k}(s_{k}) \ge \frac{\Delta t_{k}}{4(1+\Delta t_{k})}
    \left\|p_{g_{k}} \right\|^{2}
    \ge c_{w} \|p_{g_k}\| \|s_{k}\|,     \label{PLBREDST}
\end{align}
where $c_{m}$ is a positive constant, $p_{g_k} = Pg_{k} = P\nabla f(x_{k})$ and
the projection matrix $P$ is defined by equation \eqref{PROMAT}.
\end{lemma}
\proof From Lemma \ref{LEMIVBFGS}, the L-BFGS formula \eqref{LBFGSR} and the
regularization continuation method \eqref{TRCM}-\eqref{XK1}, we know that
$Ps_{k} = s_{k} \, (k = 0, \, 1, \, 2, \, \ldots)$. Furthermore, from the
L-BFGS formula and Lemma \ref{LEMHLB}, we know that the eigenvalues of
$B_{k}^{-1}$ are greater than 1/2. By combining them into equation \eqref{LOAM}
and using the symmetric Shur decomposition (p. 440, \cite{GV2013}) of
$B_{k}^{-1}$, we obtain
\begin{align}
     & m_{k}(0) - m_{k}(s_{k}) = - \frac{1 + 0.5 \Delta t_{k}}{1 + \Delta t_{k}}
     g_{k}^{T}s_{k} = - \frac{1 + 0.5 \Delta t_{k}}{1 + \Delta t_{k}} g_{k}^{T}(Ps_{k})
     = - \frac{1 + 0.5 \Delta t_{k}}{1 + \Delta t_{k}} p_{g_k}^{T}s_{k} \nonumber \\
     & \hskip 2mm
     = \frac{1 + 0.5 \Delta t_{k}}{1 + \Delta t_{k}}
     \frac{\Delta t_{k}}{1+\Delta t_{k}} p_{g_k}^{T}B_{k}^{-1}p_{g_k}
     \ge \frac{1 + 0.5 \Delta t_{k}}{1 + \Delta t_{k}}
     \frac{\Delta t_{k}}{2(1+\Delta t_{k})} \|p_{g_k}\|^{2}.
     \label{LBWPREM}
\end{align}
By using the property
$(1 + 0.5 \Delta t_{k})/(1+\Delta t_{k}) \ge (0.5+0.5\Delta t_{k})/(1+\Delta t_{k}) = 0.5$,
from equation \eqref{LBWPREM}, we have
\begin{align}
      m_{k}(0) - m_{k}(s_{k}) \ge \frac{\Delta t_{k}}{4(1+\Delta t_{k})} \|p_{g_k}\|^{2}.
     \label{LBWPPREM}
\end{align}

\vskip 2mm

From equation \eqref{BOUNDPHESS}, we have
\begin{align}
   & \|y_{k-1}\| = \|Pg(x_{k-1}) - Pg(x_{k-2})\| =
   \left\|\int_{0}^{1} P\nabla^{2}f(x_{k-2} + ts_{k-1})s_{k-1}dt\right\| \nonumber \\
   & \hskip 2mm
   = \left\|\int_{0}^{1} P\nabla^{2}f(x_{k-2} + ts_{k-1})Ps_{k-1}dt\right\| \nonumber \\
   & \hskip 2mm
   \le \int_{0}^{1} \left\|P\nabla^{2}f(x_{k-2} + ts_{k-1})P\right\| \|s_{k-1}\|dt
   \le M \|s_{k-1}\|.  \label{YKLESK}
\end{align}
From Lemma \ref{LEMHLB}, we know that the eigenvalues of $B_{k}$ are greater
than $\frac{\theta^{2} \|s_{k-1}\|^{2}}{2\|y_{k-1}\|^{2}}$. By combining it with
inequality \eqref{YKLESK}, we know that the eigenvalues of $B_{k}$ are
greater than $\theta^{2}/(2M^{2})$. Furthermore, from the symmetric Shur
decomposition (p. 440, \cite{GV2013}), we know that there exists an orthogonal
matrix $U_{k}$ such that $B_{k} = U_{k}^{T}\text{diag}
\left(\mu_{1}^{k}, \, \ldots, \, \mu_{n}^{k}\right)U_{k}$, where $\mu_{1}^{k} \ge
\mu_{2}^{k} \ge \cdots \ge \mu_{n}^{k}$ are the eigenvalues of the symmetric
matrix $B_{k}$. Thus, we obtain
\begin{align}
    & \|B_{k}s_{k}\|^{2} = \left\|(U_{k}B_{k}U_{k}^{T})U_{k}s_{k}\right\|^{2}
      = (U_{k}s_{k})^{T}\text{diag}\left((\mu_{1}^{k})^{2}, \, \ldots, \,
      (\mu_{n}^{k})^{2}\right)(U_{k}s_{k}) \nonumber \\
    & \hskip 2mm
    \ge \left(\frac{\theta^{2}}{2M^{2}}\right)^{2}  s_{k}^{T}U_{k}^{T}U_{k}s_{k}
    = \left(\frac{\theta^{2}}{2M^{2}}\right)^{2} \|s_{k}\|^{2}, \nonumber
\end{align}
which gives
\begin{align}
    \|B_{k}s_{k}\| \ge \frac{\theta^{2}}{2M^{2}} \|s_{k}\|. \label{BKSKGE}
\end{align}
By combining it with equations \eqref{TRCM} and \eqref{LBWPPREM},  we obtain
\begin{align}
     m_{k}(0) - m_{k}(s_{k}) \ge
     \frac{\Delta t_{k}}{4(1+\Delta t_{k})} \|p_{g_k}\|^{2}
     = \frac{1}{4} \|p_{g_k}\| \|B_{k}s_{k}\|
     \ge \frac{\theta^{2}}{8M^{2}}\|p_{g_k}\| \|s_{k}\|.  \label{LBWPPREDM}
\end{align}
We set $c_{w} = \theta^{2}/(8M^{2})$. Then, from equation \eqref{LBWPPREDM},
we obtain the result \eqref{PLBREDST}. \qed

\vskip 2mm

\begin{lemma} \label{LBSOAMIP}
Assume that the approximation model $m_{k}(s)$ is defined by equation \eqref{LOAM}
and $s_{k}$ is computed by the regularization continuation method
\eqref{TRCM}-\eqref{XK1}, where $B_{k} =
\left(\frac{\sigma_{0}}{\Delta t_{k}}I + P\nabla^{2} f(x_{k}) P\right)$ and
$\Delta t_{k} \le \frac{\sigma_{0}}{2M}$ in the  ill-posed phase. Then, we have
\begin{align}
    m_{k}(0) - m_{k}(s_{k}) \ge \frac{\Delta t_{k}}{4(1+\Delta t_{k})}
    \left\|p_{g_{k}} \right\|^{2}
    \ge c_{b} \|p_{g_k}\| \|s_{k}\|,     \label{PLBREDSTIP}
\end{align}
where $c_{b}$ is a positive constant, $p_{g_k} = Pg_{k} = P\nabla f(x_{k})$ and
the projection matrix $P$ is defined by equation \eqref{PROMAT}.
\end{lemma}
\proof From equations \eqref{TRCM}-\eqref{XK1} and
$B_{k} = \left(\frac{\sigma_{0}}{\Delta t_{k}} I + P\nabla^{2}f(x_{k})P\right)$, we
have
\begin{align}
     P\left(\frac{\sigma_{0}}{\Delta t_{k}}I + P\nabla^{2}f(x_{k})P\right) s_{k}
     = - \frac{\Delta t_{k}}{1+ \Delta t_{k}}P^{2}g_{k}. \label{PTRCM}
\end{align}
By substituting $P^{2} = P$ into equation \eqref{PTRCM}, we obtain
$Ps_{k} = s_{k}$. Consequently, by combining it with the property $AP = 0$,
we obtain $As_{k} = 0$,  i.e. $x_{k+1} \in S_{f}$ if $x_{k} \in S_{f}$.
By induction, we obtain $x_{k} \in S_{f} \, (k = 1, \, 2, \, \ldots)$ when
$x_{0} \in S_{f}$. Therefore, according to the assumption
$\Delta t_{k} \le \frac{\sigma_{0}}{2M}$, from equation \eqref{SPDCON}, we know
\begin{align}
     \left(\frac{\sigma_{0}}{\Delta t_{k}}I + P\nabla^{2} f(x_{k})P\right) \succ 0.
     \label{TRPSD}
\end{align}

\vskip 2mm

From equations \eqref{TRCM}, \eqref{TRPSD} and $Ps_{k} = s_{k}$, by using the
symmetric Shur decomposition (p. 440, \cite{GV2013}), we have
\begin{align}
    & - s_{k}^{T}g_{k} = - (Ps_{k})^{T}g_{k}
    = - s_{k}^{T}(Pg_{k}) = \frac{\Delta t_{k}}{1+\Delta t_{k}}
    p_{g_k}^{T}\left(\frac{\sigma_{0}}{\Delta t_{k}}I + P\nabla^{2}f(x_{k})P\right)^{-1}p_{g_k}
    \nonumber \\
    & \hskip 2mm
    \ge \frac{\Delta t_{k}}{1+\Delta t_{k}}
    \frac{1}{\sigma_{0}/\Delta t_{k}+\|P\nabla^{2}f(x_{k})P\|}\|p_{g_k}\|^{2}
    \ge  \frac{\Delta t_{k}}{1+\Delta t_{k}}
    \frac{1}{\sigma_{0}/\Delta t_{k}+M}\|p_{g_k}\|^{2}.     \label{SKINGKGE}
\end{align}
Similarly to the estimation of equation \eqref{BKSKGE}, from equation \eqref{TRCM}
and the symmetric Shur decomposition (p. 440, \cite{GV2013}), we have
\begin{align}
     \|B_{k}s_{k}\| =
     \left\|\left(\frac{\sigma_{0}}{\Delta t_{k}}I + P\nabla^{2}f(x_{k})P\right)s_{k}\right\|
     \ge \left(\frac{\sigma_{0}}{\Delta t_{k}} - M\right)\|s_{k}\|, \label{BKSKGELB}
\end{align}
where we use the property that the absolute eigenvalues of $P\nabla^{2}f(x_{k})P$
are less than $M$. From equations \eqref{TRCM} and \eqref{SKINGKGE}-\eqref{BKSKGELB},
we obtain
\begin{align}
   - s_{k}^{T}g_{k} & \ge \frac{\|B_{k}s_{k}\| \|p_{g_k}\|}{\sigma_{0}/\Delta t_{k} + M}
   \ge \frac{\sigma_{0}/\Delta t_{k} - M }{\sigma_{0}/\Delta t_{k} + M}
   \|p_{g_k}\| \|s_{k}\| \nonumber \\
  & \ge \frac{2M - M }{2M + M}
   \|p_{g_k}\| \|s_{k}\| = \frac{1}{3} \|p_{g_k}\| \|s_{k}\|,
  \label{SKGKGELB}
\end{align}
where we use the assumption $\Delta t_{k} \le \sigma_{0}/(2M)$ and the monotonically
increasing property of $\alpha(t) = (t - M)/(t +M)$ when $t > M$.

\vskip 2mm

From the approximation model \eqref{LOAM} and the estimation \eqref{SKGKGELB},
we have
\begin{align}
    & m_{k}(0) - m_{k}(s_{k}) = - \frac{1 + 0.5 \Delta t_{k}}{1+ \Delta t_{k}}
    g_{k}^{T}s_{k}
    \ge \frac{1}{3} \frac{1 + 0.5 \Delta t_{k}}{1+ \Delta t_{k}}\|p_{g_k}\| \|s_{k}\|
    \nonumber \\
    & \hskip 2mm
    =  \frac{1}{3} \frac{0.5 + 0.5 (1+\Delta t_{k})}{1+ \Delta t_{k}}\|p_{g_k}\| \|s_{k}\|
    \ge \frac{1}{6} \|p_{g_k}\| \|s_{k}\|, \label{IPMODELGE}
\end{align}
where we use the property $0.5 + 0.5 (1+\Delta t_{k}) \ge 0.5 (1+\Delta t_{k})$.
We set $c_{b} = 1/6$. Then, from equation \eqref{IPMODELGE}, we obtain the
estimation \eqref{PLBREDSTIP}. \qed

\vskip 2mm

In order to prove that $p_{g_k}$ converges to zero when $k$ tends to infinity,
we need to estimate the lower bound of time steps
$\Delta t_{k} \, (k = 1, \, 2, \, \ldots)$.

\vskip 2mm

\begin{lemma} \label{DTBOUND}
Assume that $f$ satisfies Assumption \ref{ASSFUN} and the sequence $\{x_{k}\}$ is
generated by Algorithm \ref{ALGTRCMTR}. Then, there exists a positive constant
$\delta_{\Delta t}$ such that
\begin{align}
    \Delta t_{k} \ge \gamma_{2} \delta_{\Delta t}     \label{DTGEPN}
\end{align}
holds for all $k = 1, \,  2, \, \dots$, where $\Delta t_{k}$ is adaptively adjusted
by the trust-region updating strategy \eqref{LOAM}-\eqref{ADTK1}.
\end{lemma}

\vskip 2mm

\proof From the first-order Taylor expansion, we have
\begin{align}
    f(x_{k}+ s_{k}) =  f(x_{k}) +
    \int_{0}^{1} s_{k}^{T}g (x_{k} + t s_{k}) dt.\label{FOTEFK}
\end{align}
Thus, from equations \eqref{LOAM}-\eqref{MRHOK}, \eqref{FOTEFK}, the Armijo
sufficient descent condition \eqref{AMGESKPGK} and the assumption
\eqref{BOUNDHESS}, we have
\begin{align}
      & \left|\rho_{k} - 1\right| =  \left|\frac{(f(x_{k}) - f(x_{k}+s_{k}))
       - (m_{k}(0) - m_{k}(s_{k}))}{m_{k}(0) - m_{k}(s_{k})}\right|
       \nonumber \\
      & \quad \le
      \frac{\left|\int_{0}^{1}s_{k}^{T}(g(x_{k} + t s_{k}) - g(x_{k}))dt\right|}
      {m_{k}(0) - m_{k}(s_{k})} + \frac{0.5\Delta t_{k}}{1+0.5\Delta t_{k}}
      \nonumber \\
      & \quad \le
      \frac{\int_{0}^{1}\|s_{k}\| \|g(x_{k} + t s_{k}) - g(x_{k})\|dt}
      {m_{k}(0) - m_{k}(s_{k})} + \frac{0.5\Delta t_{k}}{1+0.5\Delta t_{k}}
      \nonumber \\
      & \quad \le  \frac{0.5M\|s_{k}\|^{2}}{m_{k}(0) - m_{k}(s_{k})}
       + \frac{0.5\Delta t_{k}}{1+0.5\Delta t_{k}}.
      \label{ESTRHOK}
\end{align}

\vskip 2mm

From Lemma \ref{LBSOAM} and Lemma \ref{LBSOAMIP}, we know that there exists a
constant $\eta_{m}$ such as $\eta_{m} = \min\{c_{w}, \, c_{b}\}$  such that the
approximation model $m_{k}(0) - m_{k}(s_{k})$ satisfies the Armijo sufficient
descent condition \eqref{AMGESKPGK} when $\Delta t_{k} \le 1/(2M)$ and
$f(\cdot)$ satisfies Assumption \ref{ASSFUN}. By substituting the
sufficient descent condition  \eqref{AMGESKPGK} into equation \eqref{ESTRHOK},
we obtain
\begin{align}
     \left|\rho_{k} - 1\right| \le \frac{0.5M}{\eta_{m}}
     \frac{\|s_{k}\|}{\|p_{g_k}\|}
     + \frac{0.5\Delta t_{k}}{1+0.5\Delta t_{k}}. \label{UBROHK}
\end{align}

\vskip 2mm

When $B_{k}$ is updated by the L-BFGS formula \eqref{LBFGSR} in the
well-posed phase,  from Lemma \ref{LEMHLB}, we know that the eigenvalues of
$B_{k}^{-1}$ are less than
$\max\left\{1, \; \frac{2\|y_{k-1}\|^{2}}{\theta^{2}\|s_{k-1}\|^{2}}\right\}$.
By combining it with equations \eqref{LBFGSR} and \eqref{BOUNDPHESS}, we obtain
\begin{align}
     & \|s_{k}\| = \frac{\Delta t_{k}}{1+ \Delta t_{k}}\left\|B_{k}^{-1}p_{g_k}\right\|
     \le \frac{\Delta t_{k}}{1+ \Delta t_{k}}
     \max\left\{1, \; \frac{2\|y_{k-1}\|^{2}}{\theta^{2}\|s_{k-1}\|^{2}}\right\}
     \|p_{g_k}\| \nonumber \\
     \hskip 2mm
     & \le \frac{\Delta t_{k}}{1+ \Delta t_{k}}
     \max\left\{1, \; \frac{2\|Pg(x_{k-1}+s_{k-1}) - Pg(x_{k-1})\|^{2}}
     {\theta^{2}\|s_{k-1}\|^{2}}\right\}\|p_{g_k}\|     \nonumber \\
     \hskip 2mm
     & \le \frac{\Delta t_{k}}{1+ \Delta t_{k}}
     \max\left\{1, \; \frac{2M^{2}}{\theta^{2}}\right\}
     \|p_{g_k}\| = \frac{\Delta t_{k}}{1+ \Delta t_{k}}L_{w} \|p_{g_k}\|,
     \label{WPUBSK}
\end{align}
where $L_{w} \triangleq \max\left\{1, \; 2M^{2}/\theta^{2}\right\}$.

\vskip 2mm

When $B_{k} = \left(\frac{\sigma_{0}}{\Delta t_{k}}I + P\nabla^{2}f(x_k)P \right)$
and $\Delta t_{k} \le \sigma_{0}/(2M)$, from equations \eqref{TRCM} and
\eqref{BOUNDPHESS}, we have
\begin{align}
     & \|s_{k}\| = \frac{\Delta t_{k}}{1+\Delta t_{k}}\left\|B_{k}^{-1}p_{g_k}\right\|
     = \frac{\Delta t_{k}}{1+\Delta t_{k}}
     \left\|\left(\frac{\sigma_{0}}{\Delta t_{k}}I + P\nabla^{2}f(x_k)P \right)^{-1}p_{g_k}\right\|
     \nonumber \\
     \hskip 2mm
     & \le \frac{\Delta t_{k}}{1+\Delta t_{k}} \frac{1}{\sigma_{0}/\Delta t_{k} - M} \|p_{g_k}\|
     \le \frac{\Delta t_{k}}{1+\Delta t_{k}} \frac{1}{M} \|p_{g_k}\|.
     \label{IPBUBSK}
\end{align}
Thus, when $B_{k}$ are updated by the formula \eqref{UPDJK1} and
$\Delta t_{k} \le \sigma_{0}/(2M)$ in the ill-posed phase, from equation \eqref{IPBUBSK},
we have
\begin{align}
    \|s_k\| \le \frac{\Delta t_{k}}{1+\Delta t_{k}} \frac{1}{M} \|p_{g_k}\|.
     \label{IPUBSK}
\end{align}

\vskip 2mm

We set $L_{u} \triangleq \max \{L_{w}, \, 1/M\}$. By substituting equations
\eqref{WPUBSK} and \eqref{IPUBSK} into equation \eqref{UBROHK},
when $\Delta t_{k} \le \sigma_{0}/(2M)$, we obtain
\begin{align}
     & \left|\rho_{k} - 1\right| \le \frac{0.5ML_{u}}{\eta_{m}}
     \frac{\Delta t_{k}}{1+ \Delta t_{k}}
     + \frac{0.5\Delta t_{k}}{1+0.5\Delta t_{k}}\nonumber \\
     & \hskip 2mm
     \le \frac{0.5ML_{u}}{\eta_{m}}
     \frac{\Delta t_{k}}{1+ \Delta t_{k}}
     + \frac{0.5\Delta t_{k}}{0.5+0.5\Delta t_{k}}
     \le \frac{0.5ML_{u}+\eta_{m}}{\eta_{m}}
     \frac{\Delta t_{k}}{1+ \Delta t_{k}}. \label{UBROHKA1}
\end{align}
We set
\begin{align}
     \delta_{\Delta t} \triangleq
     \min\left\{\frac{\eta_{1}\eta_{m}}{0.5ML_{u}+\eta_{m}},
     \; \frac{\sigma_{0}}{2M}, \; \Delta t_{0}\right\}.      \label{UPBMPD}
\end{align}
Then, from equations \eqref{UBROHKA1}-\eqref{UPBMPD}, when
$\Delta t_{k} \le \delta_{\Delta t}$, it is not difficult to verify
\begin{align}
    \left|\rho_{k} - 1\right|  \le \eta_{1}.     \label{RHOLETA1}
\end{align}

\vskip 2mm

We assume that $K$ is the first index such that $\Delta t_{K} \le
\delta_{\Delta t}$ where $\delta_{\Delta t}$ is defined by equation \eqref{UPBMPD}.
Then, from equations \eqref{UPBMPD}-\eqref{RHOLETA1}, we know that
$|\rho_{K} - 1 | \le \eta_{1}$. According to the time step adjustment
formula \eqref{ADTK1}, $x_{K} + s_{K}$ will be accepted and the time step
$\Delta t_{K+1}$ will be enlarged. Consequently,
$\Delta t_{k}\ge \gamma_{2}\delta_{\Delta t}$ holds for all
$k = 1, \, 2, \ldots$. \qed

\vskip 2mm

By using the result of Lemma \ref{DTBOUND}, we prove the global convergence of
Algorithm \ref{ALGTRCMTR} for the linearly constrained optimization problem
\eqref{LEQOPT} in Theorem \ref{THEGCON}.

\vskip 2mm

\begin{theorem} \label{THEGCON}
Assume that $f$ satisfies Assumption \ref{ASSFUN} and $f(x)$ is lower bounded
when $x \in S_{f}$, where $S_{f} = \{x: \; Ax = b \}$. The sequence $\{x_{k}\}$
is generated by Algorithm \ref{ALGTRCMTR}. Then, we have
\begin{align}
  \lim_{k \to \infty} \inf \|Pg_{k}\| = 0, \label{LIMPGKZ}
\end{align}
where $g_{k} = \nabla f(x_{k})$ and the projection matrix $P$ is defined by equation
\eqref{PROMAT}.
\end{theorem}
\proof We prove the result \eqref{LIMPGKZ} by contradiction. Assume that there
exists a positive constant $\epsilon$ such that
\begin{align}
     \|Pg_{k}\| > \epsilon \label{LBPGKEPS}
\end{align}
holds for all $k = 0, \, 1, \, 2, \, \ldots$. According to Lemma  \ref{DTBOUND}
and Algorithm \ref{ALGTRCMTR}, we know that there exists an infinite subsequence
$\{x_{k_{i}}\}$ such that the trial steps
$s_{k_i} \, (i = 1, \, 2, \, \ldots)$ are accepted. Otherwise, all steps are
rejected after a given iteration index, then the time step will keep
decreasing to zero, which contradicts \eqref{DTGEPN}. Therefore, from equations
\eqref{MRHOK}, \eqref{AMGESKPGK} and \eqref{LBPGKEPS}, we have
\begin{align}
     & f(x_{0}) - \lim_{k \to \infty} f(x_{k})
     = \sum_{k = 0}^{\infty} (f(x_{k}) - f(x_{k+1}))
     \ge \sum_{i = 0}^{\infty} (f(x_{k_i}) - f(x_{k_i}+s_{k_i}))
      \nonumber \\
     & \ge \eta_{a} \sum_{i = 0}^{\infty}
    \left(m_{k_{i}}(0) - m_{k_{i}}(s_{k_{i}})\right)
    \ge \eta_{a} \eta_{m} \sum_{i = 0}^{\infty} \|Pg_{k_i}\| \|s_{k_i}\|
    \ge \eta_{a} \eta_{m} \epsilon  \sum_{i = 0}^{\infty} \|s_{k_i}\|.
    \label{LIMSUMFK}
\end{align}
Since $f(x)$ is lower bounded when $x \in S_{f}$ and the sequence $\{f(x_{k})\}$
is monotonically decreasing, we have $\lim_{k \to \infty} f(x_{k}) = f^{\ast}$.
By substituting it into equation \eqref{LIMSUMFK}, we obtain
\begin{align}
     \lim_{i \to \infty} \|s_{k_i}\| = 0. \label{SKITOZ}
\end{align}

\vskip 2mm

When $B_{k}$ is updated by the L-BFGS formula \eqref{LEMHLB} in the well-posed
phase, from Lemma \ref{LEMHLB}, we know $\|B_{k}\| \le 2$. When $B_{k}$ is
updated by the formula \eqref{UPDJK1} in the ill-posed phase, from equations
\eqref{BOUNDPHESS} and \eqref{DTGEPN}, we know that
$\|B_{k}\| \le \left(\frac{\sigma_{0}}{\gamma_{2}\delta_{\Delta t}}+M\right)$.
We set
\begin{align}
    L_{B} \triangleq \max \left\{2, \;
    \left(\frac{\sigma_{0}}{\gamma_{2}\delta_{\Delta t}}+M\right) \right\}.
    \label{UPBK}
\end{align}
By substituting equations \eqref{DTGEPN} and \eqref{UPBK} into
equation \eqref{TRCM}, we obtain
\begin{align}
     \|Pg_{k_i}\| =
     \frac{1+\Delta t_{k_i}}{\Delta t_{k_i}}\left\|B_{k_i}s_{k_i}\right\|
     = \left(1+ \frac{1}{\Delta t_{k_i}}\right)\left\|B_{k_i}s_{k_i}\right\|
     \le \left(1+ \frac{1}{\gamma_{2} \delta_{\Delta t}}\right)
     L_{B} \left\|s_{k_i}\right\|. \label{PGKLESKI}
\end{align}
By substituting equation \eqref{PGKLESKI} into equation \eqref{SKITOZ}, we obtain
\begin{align}
     \lim_{i \to \infty} \|Pg_{k_i}\| = 0, \nonumber
\end{align}
which contradicts the assumption \eqref{LIMPGKZ}. Consequently, the result
\eqref{LIMPGKZ} is true. \qed

\vskip 2mm

\section{Numerical Experiments}

In this section, we conduct some numerical experiments to test the performance
of Algorithm \ref{ALGTRCMTR} (Rcmtr). The codes are executed by a HP
notebook with the Intel quad-core CPU and 8Gb memory in the MATLAB R2020a
environment \cite{MATLAB}. The two-sided projection $P\nabla^{2}f(x)P$ of 
Algorithm \ref{ALGTRCMTR} is approximated by the difference formula 
\eqref{NUMHESS}.

\vskip 2mm

SQP \cite{FP1963,Goldfarb1970,NW1999,Wilson1963} is the traditional-representative
method for the constrained optimization problems. Ptctr is the recent continuation
method and its computational efficiency is significantly better than that of SQP 
for linearly constrained optimization problems according to the numerical results 
in \cite{LLS2021}. Therefore, we select these two typical methods as the basis 
for comparison. The implementation code of SQP is the built-in subroutine fmincon.m 
of the MATLAB2020a environment \cite{MATLAB}. The alternating direction method of 
multipliers (ADMM \cite{BPCE2011}) is an efficient method for some convex 
optimization problems and studied by many researchers in recent years. Therefore, 
we also compare Rcmtr with ADMM for some linearly constrained convex optimization 
problems. The compared ADMM subroutine \cite{BPCE2011} is downloaded from the web 
site at \url{https://web.stanford.edu/~boyd/papers/admm/}.

\vskip 2mm

We select $57$ optimization problems from references 
\cite{AD2005,LLS2021,ML2004,SB2013} as the test problems, some of which are the
unconstrained optimization problems \cite{AD2005,ML2004,SB2013} and we add the
same linear constraint $Ax = b$, where $b = 2*\text{ones}(n,1)$ and $A$
is defined as follows:
\begin{align}
  A_{1} =
   \begin{bmatrix}
        2 & 1 & 0 & \cdots & 0 & 0 & 0 \\
        1 & 2 & 1 & \cdots & 0 & 0 & 0 \\
        0 & 1 & 2 & \cdots & 0 & 0 & 0 \\
        \vdots & \vdots & \vdots & \ddots & \vdots& \vdots & \vdots \\
        0 & 0 & 0 & \cdots & 2 & 1 & 0 \\
        0 & 0 & 0 & \cdots & 1 & 2 & 1 \\
        0 & 0 & 0 & \cdots & 0 & 1 & 2
    \end{bmatrix}, \;
   A_{2} =
   \begin{bmatrix}
        1 & 1 & 1 & \cdots & 1 & 1 & 1 \\
        2 & 2 & 2 & \cdots & 2 & 2 & 2 \\
        1 & 1 & 1 & \cdots & 1 & 1 & 1 \\
        \vdots & \vdots & \vdots & {} & \vdots& \vdots & \vdots \\
        1 & 1 & 1 & \cdots & 1 & 1 & 1 \\
        2 & 2 & 2 & \cdots & 2 & 2 & 2 \\
        1 & 1 & 1 & \cdots & 1 & 1 & 1
    \end{bmatrix}, \;
    A =
    \begin{bmatrix}
    A_{1} & , & A_{2}\\
    \end{bmatrix}.
     \label{FORA12}
\end{align}
The termination conditions of the four compared methods are all set by
\begin{align}
    & \|\nabla_{x} L(x_{k}, \, \lambda_{k})\|_{\infty} \le 1.0 \times 10^{-6},
    \label{FOOPTTOL} \\
   & \|Ax_k - b \|_{\infty} \le 1.0 \times 10^{-6}, \;
    k = 1, \, 2, \, \ldots,   \label{FEATOL}
\end{align}
where the Lagrange function $L(x, \, \lambda)$ is defined by equation \eqref{LAGFUN}
and $\lambda$ is defined by equation \eqref{LAMBDA}.

\vskip 2mm

\begin{figure}[!htbp]
      \centering
        \includegraphics[width=0.80\textwidth]{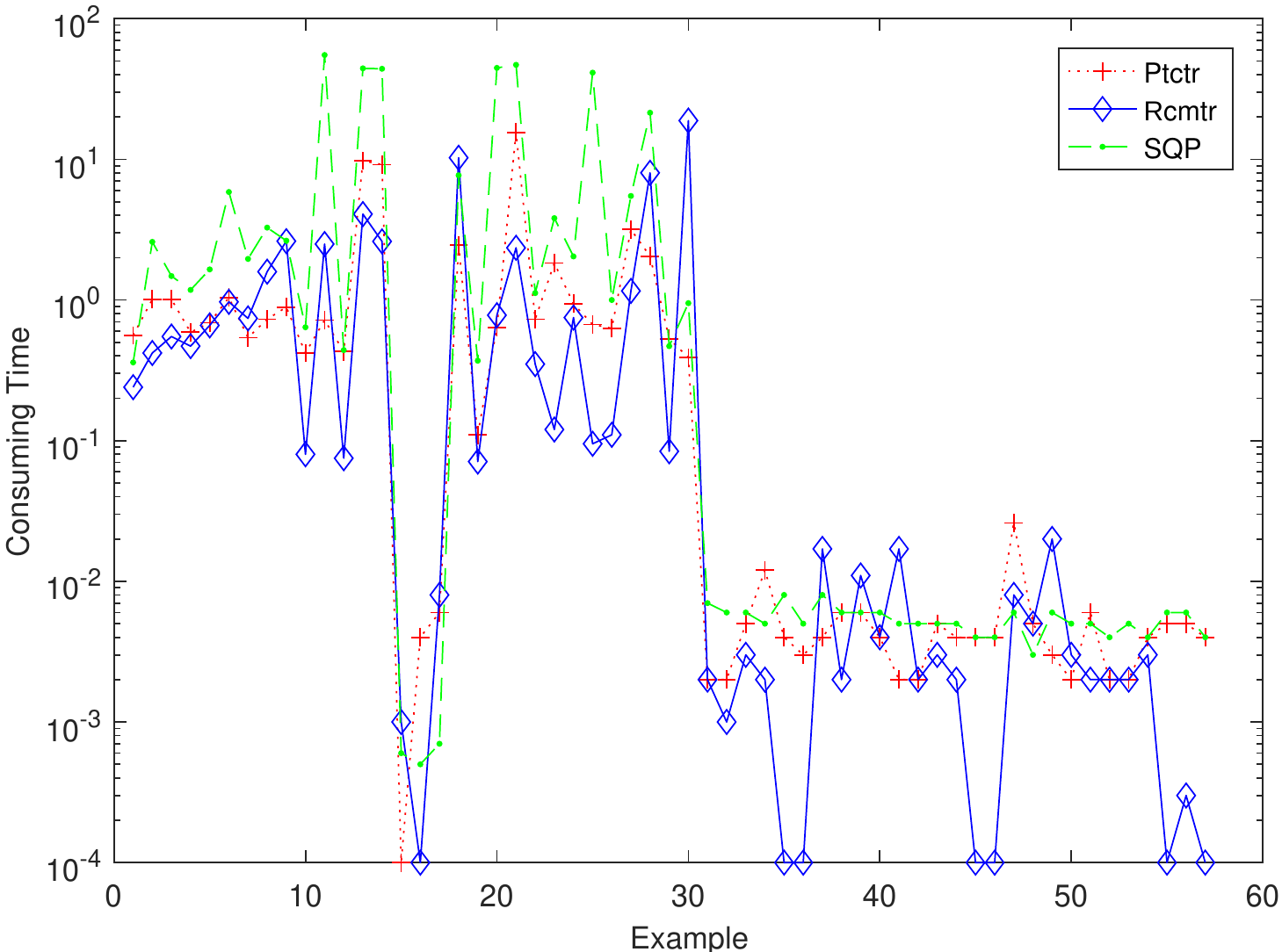}
        \caption{The computational time (s) of Ptctr, \,Rcmtr and SQP for test problems.}
        \label{fig:CNMTIM}
\end{figure}

\begin{figure}[!htbp]
      \centering
        \includegraphics[width=0.80\textwidth]{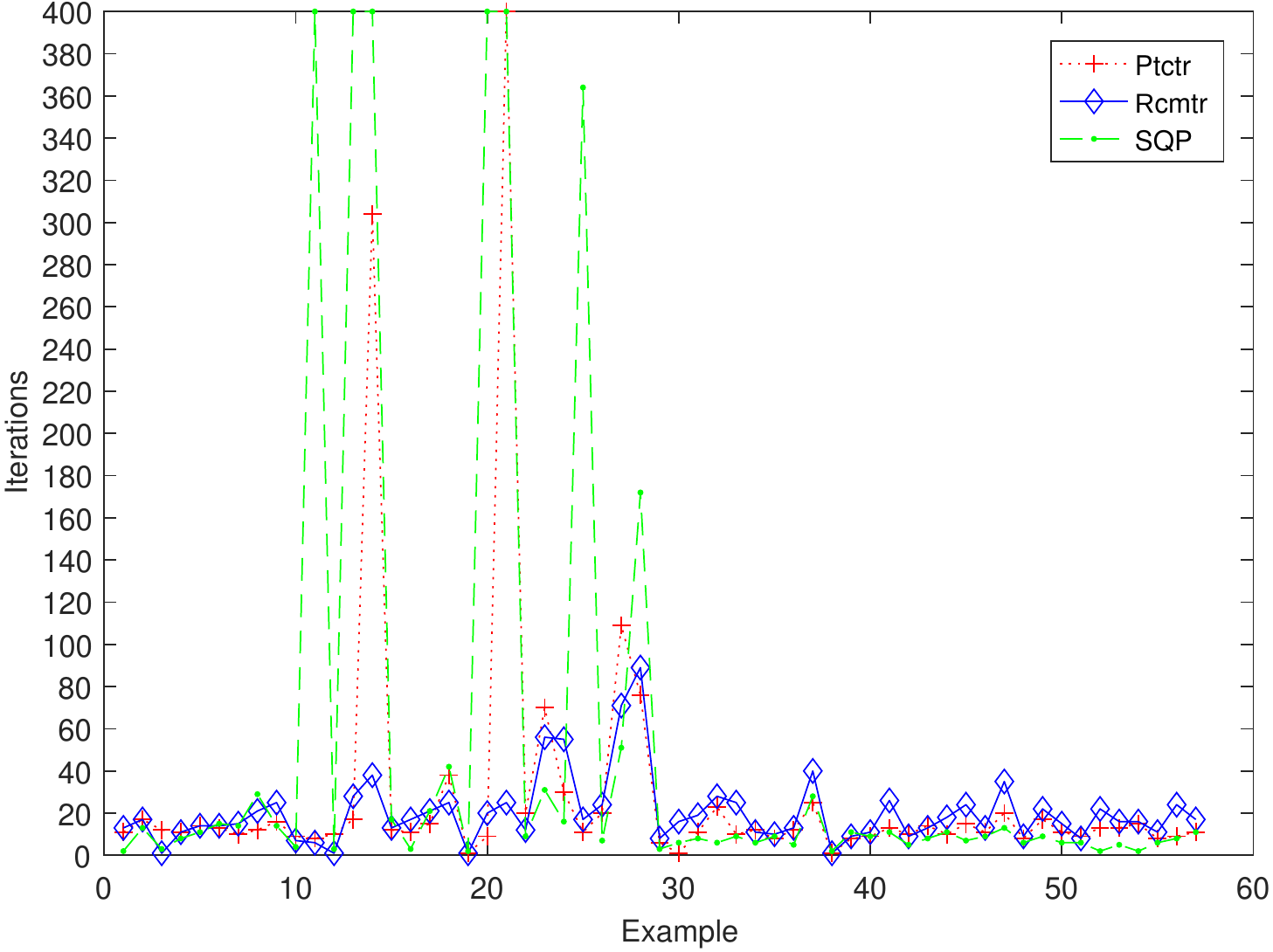}
        \caption{The number of iterations of Ptctr, \,Rcmtr and SQP for test problems.}
        \label{fig:ITESTE}
\end{figure}

We test those $57$ problems with $n = 2$ to $n \approx 1000$. The numerical
results are arranged in Tables \ref{TABCOMADMM}-\ref{TABCOMPTSC} for the convex 
problems, and Tables \ref{TABCOMPTSNL}-\ref{TABCOMPTSNS} for the non-convex 
problems. The computational time and the number of iterations of Rcmtr, Ptctr 
and SQP are illustrated in Figure \ref{fig:CNMTIM} and Figure
\ref{fig:ITESTE}, respectively. From Table \ref{TABCOMADMM} and Table 
\ref{TABCOMPTSC}, we find that Rcmtr can solve those convex optimization 
problems with linear equality constraints well. However, there are 3 convex 
problems of 17 convex test problems can not be solved by Ptctr and SQP, 
respectively. ADMM can not work well for those 17 test convex problems.

\vskip 2mm

From Table \ref{TABCOMPTSNL} and Table \ref{TABCOMPTSNS}, we find that Rcmtr can
solve those 40 non-convex linearly constrained optimization problems well
except for a particularly difficult problem (Strectched V Function \cite{SB2013}).
For this problem, Ptctr and SQP can not solve it, too. Ptctr and SQP can not 
solve two non-convex problems and five non-convex problems of 40 non-convex 
problems, respectively. Furthermore, from Tables \ref{TABCOMPTSC}-\ref{TABCOMPTSNS}
and Figure \ref{fig:CNMTIM}, we find that the computational time of Rcmtr is 
significantly less than those of Ptctr and SQP for most of test problems, 
respectively. The computational time of Rcmtr is about 1/3 of that of SQP 
(fmincon.m).

\vskip 2mm

From the numerical results, we find that Rcmtr works significantly better than the
other three methods. One of the reasons is that Rcmtr uses the
L-BFGS method \eqref{ILBFGS} as the preconditioned technique to follow their
trajectories in the well-posed phase. Consequently, Rcmtr only involves three
pairs of the inner product of two vectors and one matrix-vector product
($p_{g_{k}} = Pg_{k}$) to obtain the trial step $s_{k}$ and involves about
$(n-m)n$ flops at every iteration in the well-posed phase. However, Ptctr needs
to solve a linear system of equations with an $n \times n$ symmetric positive 
definite coefficient matrix and involves about $\frac{1}{3}n^{3}$ flops 
(p. 169, \cite{GV2013}) at every iteration. SQP needs to solve a linear 
system of equations with dimension $(m+n)$ when it solves a quadratic 
programming subproblem at every iteration (pp. 531-532, \cite{NW1999}) and 
involves about $\frac{2}{3}(m+n)^{3}$ flops (p. 116, \cite{GV2013}).

\vskip 2mm

\begin{table}[!http]
  \newcommand{\tabincell}[2]{\begin{tabular}{@{}#1@{}}#2\end{tabular}}
  \scriptsize
  \centering
  \caption{Numerical results of Rcmtr and ADMM for convex problems.} \label{TABCOMADMM}
  \begin{tabular}{|c|c|c|c|c|c|c|c|c|c|}
  \hline
  \multirow{2}{*}{Problems } & \multicolumn{2}{c|}{Rcmtr}& \multicolumn{2}{c|}{ADMM} \\ \cline{2-5} & \tabincell{c}{steps \\(time)} & \tabincell{c}{$f(x^\star)$ \\ (KKT)} & \tabincell{c}{steps \\(time)}   & \tabincell{c}{$f(x^\star)$ \\ (KKT)} \\ \hline

  \tabincell{c}{Exam. 1 Kim Problem 1 \cite{Kim2010,LLS2021} \\ (n = 1000, m = n/2)}   & \tabincell{c}{13 \\ (0.24)} &\tabincell{c}{7.27e+03 \\ (3.44e-07)} & \tabincell{c}{ \textcolor{red}{3} \\ \textcolor{red}{(0.04)}} &\tabincell{c}{\textcolor{red}{2.20e+04} \\ \textcolor{red}{(40.00)} \\ \textcolor{red}{(failed)}} \\ \hline

  \tabincell{c}{Exam. 2 LLS Problem 1 \cite{LLS2021} \\ (n = 1200, m = n/3)}   & \tabincell{c}{17 \\ (0.42)} &\tabincell{c}{1.44e+03 \\ (9.37e-07)} & \tabincell{c}{ \textcolor{red}{21} \\ \textcolor{red}{(0.07)}} &\tabincell{c}{\textcolor{red}{2.73e+03} \\ \textcolor{red}{(4.00)} \\ \textcolor{red}{(failed)}} \\ \hline

  \tabincell{c}{Exam. 3 Obsborne Problem 1 \cite{LLS2021,Osborne2016} \\ (n = 1200, m = 2/3n)}  & \tabincell{c}{1 \\ (0.55)} &\tabincell{c}{7.15e+02 \\ (1.27e-15)} & \tabincell{c}{ \textcolor{red}{60} \\ \textcolor{red}{(0.18)}} &\tabincell{c}{\textcolor{red}{8.48e+02} \\ \textcolor{red}{(2.80)} \\ \textcolor{red}{(failed)}} \\ \hline

  \tabincell{c}{Exam. 4 Mak Problem \cite{LLS2021,MAK2019} \\ (n = 1000, m = n/2)}   & \tabincell{c}{11 \\ (0.47)} &\tabincell{c}{97.96 \\ (7.74e-07)} & \tabincell{c}{ \textcolor{red}{4} \\ \textcolor{red}{(0.05)}} &\tabincell{c}{\textcolor{red}{1.32e+02} \\ \textcolor{red}{(1.00)} \\ \textcolor{red}{(failed)}} \\ \hline

  \tabincell{c}{Exam. 5 LLS Problem 2 \cite{LLS2021} \\ (n = 1000, m = n/2)}   & \tabincell{c}{14 \\ (0.66)} &\tabincell{c}{82.43 \\ (7.54e-08)} & \tabincell{c}{\textcolor{red}{12} \\ \textcolor{red}{(0.04})} &\tabincell{c}{\textcolor{red}{8.00e+03} \\ \textcolor{red}{(32.00)} \\ \textcolor{red}{(failed)}} \\ \hline

  \tabincell{c}{Exam. 6 Osborne Problem 2 \cite{LLS2021,Osborne2016} \\ (n = 1200, m = n/2)}   & \tabincell{c}{14 \\ (0.97)} &\tabincell{c}{5.14e+02 \\ (8.75e-07)} & \tabincell{c}{\textcolor{red}{60} \\ \textcolor{red}{(0.20})} &\tabincell{c}{\textcolor{red}{7.86e+02} \\ \textcolor{red}{(2.80)} \\ \textcolor{red}{(failed)}} \\ \hline

  \tabincell{c}{Exam. 7 Carlberg Problem \cite{Carlberg2009,LLS2021} \\ (n = 1000,  m = n/2)}   & \tabincell{c}{15 \\ (0.74)} &\tabincell{c}{1.19e+04 \\ (1.66e-06)} & \tabincell{c}{\textcolor{red}{3} \\ \textcolor{red}{(0.04})} &\tabincell{c}{\textcolor{red}{1.40e+04} \\ \textcolor{red}{(32.00)} \\ \textcolor{red}{(failed)}} \\ \hline

  \tabincell{c}{Exam. 8 Kim Problem 2 \cite{Kim2010,LLS2021} \\ (n = 1000, m = n/2)}   & \tabincell{c}{21 \\ (1.59)} &\tabincell{c}{4.22e+04 \\ (1.43e-06)} & \tabincell{c}{\textcolor{red}{3} \\ \textcolor{red}{(0.33})} &\tabincell{c}{\textcolor{red}{3.28e+05} \\ \textcolor{red}{(1.92e+03)} \\ \textcolor{red}{(failed)}} \\ \hline

  \tabincell{c}{Exam. 9 Yamashita Problem \cite{LLS2021,Yamashita1980} \\ (n = 1200, m = n/3)}   & \tabincell{c}{25 \\ (2.62)} &\tabincell{c}{0.50 \\ (3.67e-07)} & \tabincell{c}{\textcolor{red}{16} \\ \textcolor{red}{(0.06})} &\tabincell{c}{\textcolor{red}{25.10} \\ \textcolor{red}{(0.50)} \\ \textcolor{red}{(failed)}} \\ \hline

  \tabincell{c}{Exam. 10 Quartic With Noise \\ Function \cite{AD2005}(n = 1000,  m = n/2)}  & \tabincell{c}{7 \\ (0.08)} &\tabincell{c}{1.01e+02 \\ (2.69e-07)} & \tabincell{c}{\textcolor{red}{400} \\ \textcolor{red}{(0.40})} &\tabincell{c}{\textcolor{red}{1.01e+02} \\ \textcolor{red}{(3.98)} \\ \textcolor{red}{(failed)}} \\ \hline

  \tabincell{c}{Exam. 11 Rotated Hyper Ellopsoid \\ Function \cite{SB2013}(n = 1000,  m = n/2)}  & \tabincell{c}{6 \\ (2.50)} &\tabincell{c}{1.25e+05 \\ (8.30e-06)} & \tabincell{c}{\textcolor{red}{400} \\ \textcolor{red}{(1.04})} &\tabincell{c}{\textcolor{red}{1.26e+05} \\ \textcolor{red}{(2.00e+05)} \\ \textcolor{red}{(failed)}} \\ \hline

  \tabincell{c}{Exam. 12 Sphere Function \cite{SB2013} \\ (n = 1000,  m = n/2)}  & \tabincell{c}{1 \\ (0.08)} &\tabincell{c}{1.67e+02 \\ (3.13e-15)} & \tabincell{c}{\textcolor{red}{400} \\ \textcolor{red}{(0.27})} &\tabincell{c}{\textcolor{red}{1.67e+02} \\ \textcolor{red}{(2.00)} \\ \textcolor{red}{(failed)}} \\ \hline

  \tabincell{c}{Exam. 13 Sum Squares Function \cite{SB2013} \\ (n = 1000,  m = n/2)}  & \tabincell{c}{28 \\ (4.08)} &\tabincell{c}{4.08e+04 \\ (1.58e-06)} & \tabincell{c}{\textcolor{red}{400} \\ \textcolor{red}{(0.32})} &\tabincell{c}{\textcolor{red}{4.16e+04} \\ \textcolor{red}{(9.98e+02)} \\ \textcolor{red}{(failed)}} \\ \hline

  \tabincell{c}{Exam. 14 Trid Function \cite{SB2013} \\ (n = 1000,  m = n/2)}  & \tabincell{c}{38 \\ (2.61)} &\tabincell{c}{5.82e+02 \\ (5.36e-07)} & \tabincell{c}{\textcolor{red}{400} \\ \textcolor{red}{(0.36})} &\tabincell{c}{\textcolor{red}{5.85e+02} \\ \textcolor{red}{(3.99)} \\ \textcolor{red}{(failed)}} \\ \hline

  \tabincell{c}{Exam. 15 Booth Function \cite{SB2013} \\ (n = 2,  m = n/2)}  & \tabincell{c}{13 \\ (1.00e-03)} &\tabincell{c}{9.00 \\ (1.98e-07)} & \tabincell{c}{\textcolor{red}{18} \\ \textcolor{red}{(1.00e-03})} &\tabincell{c}{\textcolor{red}{45.00} \\ \textcolor{red}{(30.00)} \\ \textcolor{red}{(failed)}} \\ \hline

  \tabincell{c}{Exam. 16 Matyas Function \cite{SB2013} \\ (n = 2,  m = n/2)}  & \tabincell{c}{17 \\ (1.00e-04)} &\tabincell{c}{0.18 \\ (4.44e-07)} & \tabincell{c}{\textcolor{red}{18} \\ \textcolor{red}{(2.00e-03})} &\tabincell{c}{\textcolor{red}{2.60} \\ \textcolor{red}{(5.20)} \\ \textcolor{red}{(failed)}} \\ \hline

  \tabincell{c}{Exam. 17 Zakharov Function \cite{SB2013} \\ (n = 10,  m = n/2)}  & \tabincell{c}{21 \\ (8.00e-03)} &\tabincell{c}{7.31 \\ (1.65e-07)} & \tabincell{c}{\textcolor{red}{21} \\ \textcolor{red}{(1.00e-03})} &\tabincell{c}{\textcolor{red}{4.33e+02} \\ \textcolor{red}{(1.87e+03)} \\ \textcolor{red}{(failed)}} \\ \hline

\end{tabular}
\end{table}

\vskip 2mm

\begin{table}[!http]
  \newcommand{\tabincell}[2]{\begin{tabular}{@{}#1@{}}#2\end{tabular}}
  \scriptsize
  \centering
  \caption{Numerical results of Ptctr, Rcmtr and SQP for convex problems.} \label{TABCOMPTSC}
  \begin{tabular}{|c|c|c|c|c|c|c|c|c|c|}
  \hline
  \multirow{2}{*}{Problems } & \multicolumn{2}{c|}{Ptctr} & \multicolumn{2}{c|}{Rcmtr} & \multicolumn{2}{c|}{SQP}  \\ \cline{2-7}
                        & \tabincell{c}{steps \\(time)} & \tabincell{c}{$f(x^\star)$ \\ (KKT)} & \tabincell{c}{steps \\(time)}   & \tabincell{c}{$f(x^\star)$ \\ (KKT)}  & \tabincell{c}{steps \\(time)} & \tabincell{c}{$f(x^\star)$ \\ (KKT)} \\ \hline

  \tabincell{c}{Exam. 1 Kim \\ Problem 1 \cite{Kim2010,LLS2021} \\ (n = 1000, m = n/2)}   & \tabincell{c}{11 \\ (0.56)} &\tabincell{c}{7.27e+03 \\ (5.79e-08)} & \tabincell{c}{13 \\ (0.24)} &\tabincell{c}{7.27e+03 \\ (3.44e-07)} & \tabincell{c}{2 \\ (0.36)} &\tabincell{c}{7.27e+03 \\ (8.30e-13)} \\ \hline

  \tabincell{c}{Exam. 2 LLS \\ Problem 1 \cite{LLS2021} \\ (n = 1200, m = n/3)}   & \tabincell{c}{17 \\ (1.01)} &\tabincell{c}{1.44e+03 \\ (7.36e-07)} & \tabincell{c}{17 \\ (0.42)} &\tabincell{c}{1.44e+03 \\ (9.37e-07)} & \tabincell{c}{13 \\ (2.59)} &\tabincell{c}{1.44e+03 \\ (3.42e-07)} \\ \hline

  \tabincell{c}{Exam. 3 Obsborne \\ Problem 1 \cite{LLS2021,Osborne2016} \\ (n = 1200, m = 2/3n)}  & \tabincell{c}{12 \\ (1.01)} &\tabincell{c}{7.15e+02 \\ (2.30e-07)} & \tabincell{c}{1 \\ (0.55)} &\tabincell{c}{7.15e+02 \\ (1.27e-15)} & \tabincell{c}{3 \\ (1.48)} &\tabincell{c}{7.14e+02 \\ (2.22e-15)} \\ \hline

  \tabincell{c}{Exam. 4 Mak \\ Problem \cite{LLS2021,MAK2019} \\ (n = 1000, m = n/2)}   & \tabincell{c}{11 \\ (0.59)} &\tabincell{c}{97.96 \\ (3.50e-07)} & \tabincell{c}{11 \\ (0.47)} &\tabincell{c}{97.96 \\ (7.74e-07)} & \tabincell{c}{8 \\ (1.18)} &\tabincell{c}{97.96 \\ (1.34e-10)} \\ \hline

  \tabincell{c}{Exam. 5 LLS \\ Problem 2 \cite{LLS2021} \\ (n = 1000, m = n/2)}   & \tabincell{c}{14 \\ (0.69)} &\tabincell{c}{82.43 \\ (8.79e-08)} & \tabincell{c}{14 \\ (0.66)} &\tabincell{c}{82.43 \\ (7.54e-08)} & \tabincell{c}{11 \\ (1.65)} &\tabincell{c}{82.43 \\ (1.78e-09)} \\ \hline

  \tabincell{c}{Exam. 6 Osborne \\ Problem 2 \cite{LLS2021,Osborne2016} \\ (n = 1200, m = n/2)}   & \tabincell{c}{13 \\ (1.04)} &\tabincell{c}{5.14e+02 \\ (1.79e-07)} & \tabincell{c}{14 \\ (0.97)} &\tabincell{c}{5.14e+02 \\ (8.75e-07)} & \tabincell{c}{15 \\ (5.86)} &\tabincell{c}{5.14e+02 \\ (1.75e-06)} \\ \hline

  \tabincell{c}{Exam. 7 Carlberg \\ Problem \cite{Carlberg2009,LLS2021} \\ (n = 1000, m = n/2)}   & \tabincell{c}{10 \\ (0.54)} &\tabincell{c}{1.19e+04 \\ (1.23e-07)} & \tabincell{c}{15 \\ (0.74)} &\tabincell{c}{1.19e+04 \\ (1.66e-06)} & \tabincell{c}{14 \\ (1.96)} &\tabincell{c}{1.19e+04 \\ (1.13e-05)} \\ \hline

  \tabincell{c}{Exam. 8 Kim \\ Problem 2 \cite{Kim2010,LLS2021} \\ (n = 1000, m = n/2)}   & \tabincell{c}{12 \\ (0.73)} &\tabincell{c}{4.22e+04 \\ (6.14e-06)} & \tabincell{c}{21 \\ (1.59)} &\tabincell{c}{4.22e+04 \\ (1.43e-06)} & \tabincell{c}{29 \\ (3.27)} &\tabincell{c}{4.22e+04 \\ (3.05e-06)} \\ \hline

  \tabincell{c}{Exam. 9 Yamashita \\ Problem \cite{LLS2021,Yamashita1980} \\ (n = 1200, m = n/3)}   & \tabincell{c}{16 \\ (0.89)} &\tabincell{c}{0.50 \\ (4.39e-07)} & \tabincell{c}{25 \\ (2.62)} &\tabincell{c}{0.50 \\ (3.67e-07)} & \tabincell{c}{14 \\ (2.64)} &\tabincell{c}{0.50 \\ (1.01e-07)} \\ \hline

  \tabincell{c}{Exam. 10 Quartic With\\ Noise  Function \cite{AD2005}\\(n = 1000,  m = n/2)}  & \tabincell{c}{9 \\ (0.42)} &\tabincell{c}{1.01e+02 \\ (3.14e-07)} & \tabincell{c}{7 \\ (0.08)} &\tabincell{c}{1.01e+02 \\ (2.69e-07)} & \tabincell{c}{4 \\ (0.64)} &\tabincell{c}{1.01e+02 \\ (1.25e-09)} \\ \hline

  \tabincell{c}{Exam. 11 Rotated Hyper \\Ellopsoid Function \cite{SB2013}\\(n = 1000,  m = n/2)}  & \tabincell{c}{\textcolor{red}{8} \\ \textcolor{red}{(0.72)}} &\tabincell{c}{\textcolor{red}{1.25e+05} \\ \textcolor{red}{(2.08e-04)} \\ \textcolor{red}{(failed)}} & \tabincell{c}{6 \\ (2.50)} &\tabincell{c}{1.25e+05 \\ (8.30e-06)} & \tabincell{c}{\textcolor{red}{400} \\ \textcolor{red}{(55.18)}} &\tabincell{c}{\textcolor{red}{1.46e+05} \\ \textcolor{red}{(3.22e+02)} \\ \textcolor{red}{(failed)}} \\ \hline

  \tabincell{c}{Exam. 12 Sphere \\Function \cite{SB2013} \\ (n = 1000,  m = n/2)} & \tabincell{c}{10 \\ (0.43)} &\tabincell{c}{1.67e+02 \\ (1.11e-07)} & \tabincell{c}{1 \\ (7.50e-02)} &\tabincell{c}{1.67e+02 \\ (3.13e-15)} & \tabincell{c}{3 \\ (0.44)} &\tabincell{c}{1.67e+02 \\ (7.67e-10)} \\ \hline

  \tabincell{c}{Exam. 13 Sum Squares \\Function \cite{SB2013}\\ (n = 1000,  m = n/2)} & \tabincell{c}{\textcolor{red}{17} \\ \textcolor{red}{(9.77)}} &\tabincell{c}{\textcolor{red}{4.08e+04} \\ \textcolor{red}{(1.85e-04)} \\ \textcolor{red}{(failed)}} & \tabincell{c}{28 \\ (4.08)} &\tabincell{c}{4.08e+04 \\ (1.58e-06)} & \tabincell{c}{\textcolor{red}{400} \\ \textcolor{red}{(44.36)}} &\tabincell{c}{\textcolor{red}{4.10e+04} \\ \textcolor{red}{(1.01e+02)} \\ \textcolor{red}{(failed)}} \\ \hline

  \tabincell{c}{Exam. 14 Trid \\Function \cite{SB2013} \\ (n = 1000,  m = n/2)}  & \tabincell{c}{\textcolor{red}{304} \\ \textcolor{red}{(9.18)}} &\tabincell{c}{\textcolor{red}{5.82e+02} \\ \textcolor{red}{(8.34e-04)} \\ \textcolor{red}{(failed)}} & \tabincell{c}{38 \\ (2.61)} &\tabincell{c}{5.82e+02 \\ (5.36e-07)} & \tabincell{c}{\textcolor{red}{400} \\ \textcolor{red}{(44.05)}} &\tabincell{c}{\textcolor{red}{5.82e+02} \\ \textcolor{red}{(1.56e-04)} \\ \textcolor{red}{(failed)}} \\ \hline

  \tabincell{c}{Exam. 15 Booth \\Function \cite{SB2013} \\ (n = 2,  m = n/2)}  & \tabincell{c}{12 \\ (1.00e-04)} &\tabincell{c}{9.00 \\ (1.74e-07)} & \tabincell{c}{13 \\ (1.00e-03)} &\tabincell{c}{9.00 \\ (1.98e-07)} & \tabincell{c}{17 \\ (6.00e-03)} &\tabincell{c}{9.00 \\ (3.55e-15)} \\ \hline

  \tabincell{c}{Exam. 16 Matyas \\Function\cite{SB2013} \\ (n = 2,  m = n/2)}  & \tabincell{c}{11 \\ (4.00e-03)} &\tabincell{c}{0.18 \\ (1.87e-08)} & \tabincell{c}{17 \\ (1.00e-04)} &\tabincell{c}{0.18 \\ (4.44e-07)} & \tabincell{c}{3 \\ (5.00e-03)} &\tabincell{c}{0.18 \\ (1.67e-16)} \\ \hline

  \tabincell{c}{Exam. 17 Zakharov \\Function\cite{SB2013} \\ (n = 10,  m = n/2)} & \tabincell{c}{15 \\ (6.00e-03)} &\tabincell{c}{7.31 \\ (2.93e-08)} & \tabincell{c}{21 \\ (8.00e-03)} &\tabincell{c}{7.31 \\ (1.65e-07)} & \tabincell{c}{21 \\ (7.00e-03)} &\tabincell{c}{7.31 \\ (8.50e-06)} \\ \hline

\end{tabular}
\end{table}

\vskip 2mm

\begin{table}[!http]
    \newcommand{\tabincell}[2]{\begin{tabular}{@{}#1@{}}#2\end{tabular}}
    \scriptsize
    \centering
    \caption{Numerical results of Ptctr, Rcmtr, SQP for large-scale nonconvex problems.}
    \label{TABCOMPTSNL}
    \begin{tabular}{|c|c|c|c|c|c|c|c|c|c|}
      \hline
        \multirow{2}{*}{Problems } & \multicolumn{2}{c|}{Ptctr} & \multicolumn{2}{c|}{Rcmtr} & \multicolumn{2}{c|}{SQP}  \\ \cline{2-7} & \tabincell{c}{steps \\(time)} & \tabincell{c}{$f(x^\star)$ \\ (KKT)} & \tabincell{c}{steps \\(time)}
        & \tabincell{c}{$f(x^\star)$ \\ (KKT)}  & \tabincell{c}{steps \\(time)} & \tabincell{c}{$f(x^\star)$ \\ (KKT)}
        \\ \hline

      \tabincell{c}{Exam. 18 LLS \\ Problem 3 \cite{LLS2021} \\ (n = 1000, m = n/2)}   & \tabincell{c}{38 \\ (2.45)}
      &\tabincell{c}{1.96e+02 \\ (1.17e-05)} & \tabincell{c}{25 \\ (10.27)} &\tabincell{c}{-3.03e+03 \\ (4.86e-07)}
      & \tabincell{c}{42 \\ (7.70)} &\tabincell{c}{1.88e+02 \\ (7.97e-06)} \\ \hline

       \tabincell{c}{Exam. 19 Ackly \\ Function \cite{SB2013} \\ (n = 1000, m = n/2)}   & \tabincell{c}{1 \\ (0.11)} &\tabincell{c}{2.64 \\ (1.87e-07)} & \tabincell{c}{1 \\ (7.10e-02)} &\tabincell{c}{2.64 \\ (7.50e-07)}
       & \tabincell{c}{2 \\ (0.37)} &\tabincell{c}{2.42 \\ (1.94e-07)} \\ \hline

       \tabincell{c}{Exam. 20 Rosenbrock \\ Function \cite{SB2013} \\ (n = 1000, m = n/2)}   & \tabincell{c}{9 \\ (0.64)} &\tabincell{c}{9.26e+03 \\ (9.03e-06)} & \tabincell{c}{20 \\ (0.78)} &\tabincell{c}{9.26e+03 \\ (2.15e-06)} & \tabincell{c}{\textcolor{red}{400} \\ \textcolor{red}{(44.68)}} &\tabincell{c}{\textcolor{red}{9.26e+03} \\ \textcolor{red}{(5.00e-03)} \\ \textcolor{red}{(failed)}} \\ \hline

       \tabincell{c}{Exam. 21 Dixon-Price \\ Function \cite{SB2013} \\ (n = 1000, m = n/2)}
       & \tabincell{c}{\textcolor{red}{400} \\ \textcolor{red}{(15.54)}} &\tabincell{c}{\textcolor{red}{8.97e+04} \\ \textcolor{red}{(2.42e-02)} \\ \textcolor{red}{(failed)}} & \tabincell{c}{25 \\ (2.35)}
       &\tabincell{c}{9.00e+04 \\ (1.74e-09)} & \tabincell{c}{\textcolor{red}{400} \\ \textcolor{red}{(46.97)}} &\tabincell{c}{\textcolor{red}{8.24e+06} \\ \textcolor{red}{(1.28e+05)} \\ \textcolor{red}{(failed)}} \\ \hline

       \tabincell{c}{Exam. 22 Griewank \\ Function \cite{SB2013} \\ (n = 1000, m = n/2)}   & \tabincell{c}{20 \\ (0.73)}
       &\tabincell{c}{0.86 \\ (4.81e-07)} & \tabincell{c}{12 \\ (0.35)} &\tabincell{c}{0.86 \\ (4.40e-08)}
       & \tabincell{c}{9 \\ (1.12)} &\tabincell{c}{0.86 \\ (1.07e-10)} \\ \hline

       \tabincell{c}{Exam. 23 Levy \\ Function \cite{SB2013} \\ (n = 1000, m = n/2)}   & \tabincell{c}{70 \\ (1.83)}
       &\tabincell{c}{71.06 \\ (2.36e-08)} & \tabincell{c}{56 \\ (0.12)} &\tabincell{c}{71.06 \\ (8.25e-07)}
       & \tabincell{c}{31 \\ (3.82)} &\tabincell{c}{71.06 \\ (1.11e-07)} \\ \hline

       \tabincell{c}{Exam. 24 Molecular \\ Energy Function \cite{ML2004} \\ (n = 1000, m = n/2)}
       & \tabincell{c}{30 \\ (0.94)} &\tabincell{c}{4.69e+02 \\ (4.38e-07)} & \tabincell{c}{55 \\ (0.75)} &\tabincell{c}{4.69e+02 \\ (8.66e-07)} & \tabincell{c}{16 \\ (2.04)} &\tabincell{c}{4.69e+02 \\ (1.71e-06)} \\ \hline

      \tabincell{c}{Exam. 25 Powell \\ Function \cite{SB2013} \\ (n = 1000, m = n/2)}   & \tabincell{c}{11 \\ (0.67)} &\tabincell{c}{4.26e+03 \\ (1.77e-06)} & \tabincell{c}{17 \\ (9.50e-02)} &\tabincell{c}{4.26e+03 \\ (4.52e-07)}
      & \tabincell{c}{\textcolor{red}{364} \\ \textcolor{red}{(41.34)}} &\tabincell{c}{\textcolor{red}{4.26e+03} \\ \textcolor{red}{(1.38e-04)} \\ \textcolor{red}{(failed)}} \\ \hline

      \tabincell{c}{Exam. 26 Rastrigin \\ Function \cite{SB2013} \\ (n = 1000, m = n/2)}   & \tabincell{c}{20 \\ (0.63)}
      &\tabincell{c}{2.93e+03 \\ (6.42e-07)} & \tabincell{c}{24 \\ (0.11)} &\tabincell{c}{2.93e+03 \\ (1.56e-06)} & \tabincell{c}{7 \\ (1.00)} &\tabincell{c}{4.44e+03 \\ (1.13e-06)} \\ \hline

      \tabincell{c}{Exam. 27 Schwefel \\ Function \cite{SB2013} \\ (n = 1000, m = n/2)} & \tabincell{c}{109 \\ (3.17)}
      &\tabincell{c}{4.19e+05 \\ (3.74e-07)} & \tabincell{c}{71 \\ (1.16)} &\tabincell{c}{4.19e+05 \\ (4.16e-06)} & \tabincell{c}{51 \\ (5.49)} &\tabincell{c}{4.02e+05 \\ (2.31e-06)} \\ \hline

      \tabincell{c}{Exam. 28 Styblinski \\ Tang Function \cite{SB2013} \\ (n = 1000, m = n/2)} & \tabincell{c}{76 \\ (2.05)}
      &\tabincell{c}{-9.61e+03 \\ (1.12e-05)} & \tabincell{c}{89 \\ (8.02)} &\tabincell{c}{-9.61e+03 \\ (4.57e-06)}
      &\tabincell{c}{\textcolor{red}{172} \\ \textcolor{red}{(21.46)}} &\tabincell{c}{\textcolor{red}{-2.56e+04} \\
      \textcolor{red}{(7.03e-04)} \\ \textcolor{red}{(failed)}} \\ \hline

      \tabincell{c}{Exam. 29 Shubert \\ Function \cite{SB2013} \\ (n = 1000, m = n/2)}   & \tabincell{c}{6 \\ (0.53)}
      &\tabincell{c}{2.65e+03 \\ (1.43e-06)} & \tabincell{c}{8 \\ (8.40e-02)} &\tabincell{c}{2.65e+03 \\ (8.92e-07)}
      &\tabincell{c}{3 \\ (0.47)} &\tabincell{c}{2.65e+03 \\ (1.42e-05)} \\ \hline

      \tabincell{c}{Exam. 30 Strectched \\ V Function \cite{SB2013} \\ (n = 1000, m = n/2)}
      & \tabincell{c}{\textcolor{red}{1} \\ \textcolor{red}{(0.39)}} &\tabincell{c}{\textcolor{red}{3.10e-03} \\ \textcolor{red}{(1.08e+05)} \\ \textcolor{red}{(failed)}} & \tabincell{c}{\textcolor{red}{16} \\ \textcolor{red}{(18.86)}}
      &\tabincell{c}{\textcolor{red}{2.30e+02} \\ \textcolor{red}{(2.41e-03)} \\ \textcolor{red}{(failed)}}
      &\tabincell{c}{\textcolor{red}{6} \\ \textcolor{red}{(0.95)}} &\tabincell{c}{\textcolor{red}{1.28} \\ \textcolor{red}{(34.15)} \\ \textcolor{red}{(failed)}} \\ \hline

  \end{tabular}
\end{table}

\vskip 2mm

\begin{table}[!http]
  \newcommand{\tabincell}[2]{\begin{tabular}{@{}#1@{}}#2\end{tabular}}
  \scriptsize
  \centering
  \caption{Numerical results of Ptctr, Rcmtr, SQP for small-scale nonconvex problems.} \label{TABCOMPTSNS}
  \begin{tabular}{|c|c|c|c|c|c|c|c|c|c|}
  \hline
  \multirow{2}{*}{Problems } & \multicolumn{2}{c|}{Ptctr} & \multicolumn{2}{c|}{Rcmtr} & \multicolumn{2}{c|}{SQP}  \\ \cline{2-7}
  & \tabincell{c}{steps \\(time)} & \tabincell{c}{$f(x^\star)$ \\ (KKT)} & \tabincell{c}{steps \\(time)}   & \tabincell{c}{$f(x^\star)$ \\ (KKT)}  & \tabincell{c}{steps \\(time)} & \tabincell{c}{$f(x^\star)$ \\ (KKT)} \\ \hline

  \tabincell{c}{Exam. 31 Beale\\ Function \cite{SB2013} \\ (n = 2, m = n/2)}   & \tabincell{c}{11 \\ (2.00e-03)} &\tabincell{c}{3.35 \\ (2.23e-08)} & \tabincell{c}{19 \\ (2.00e-03)} &\tabincell{c}{3.35 \\ (8.92e-07)} & \tabincell{c}{8 \\ (7.00e-03)} &\tabincell{c}{3.35 \\ (1.50e-06)} \\ \hline

  \tabincell{c}{Exam. 32 Branin\\ Function \cite{SB2013} \\ (n = 2, m = n/2)}   & \tabincell{c}{23 \\ (2.00e-03)} &\tabincell{c}{15.90 \\ (1.51e-08)} & \tabincell{c}{28 \\ (1.00e-03)} &\tabincell{c}{15.90 \\ (3.53e-07)} & \tabincell{c}{6 \\ (6.00e-03)} &\tabincell{c}{34.37 \\ (4.55e-06)} \\ \hline

  \tabincell{c}{Exam. 33 Eason\\ Function \cite{SB2013} \\ (n = 2, m = n/2)}   & \tabincell{c}{10 \\ (5.00e-03)} &\tabincell{c}{-4.19e-06 \\ (4.08e-06)} & \tabincell{c}{25 \\ (3.00e-03)} &\tabincell{c}{-4.19e-06 \\ (1.06e-07)} & \tabincell{c}{9 \\ (6.00e-03)} &\tabincell{c}{-4.55e-06 \\ (5.80e-07)} \\ \hline

  \tabincell{c}{Exam. 34 Hosaki\\ Function \cite{AD2005} \\ (n = 2, m = n/2)}   & \tabincell{c}{12 \\ (1.20e-02)} &\tabincell{c}{-0.56 \\ (4.21e-07)} & \tabincell{c}{11 \\ (2.00e-03)} &\tabincell{c}{-0.56 \\ (3.24e-08)} & \tabincell{c}{6 \\ (5.00e-03)} &\tabincell{c}{-0.56 \\ (9.77e-07)} \\ \hline

  \tabincell{c}{Exam. 35 Levy\\ Function N. 13 \cite{SB2013} \\ (n = 2, m = n/2)}   & \tabincell{c}{8 \\ (4.00e-03)} &\tabincell{c}{0.63 \\ (9.51e-07)} & \tabincell{c}{10 \\ (1.00e-04)} &\tabincell{c}{0.63 \\ (1.61e-07)} & \tabincell{c}{9 \\ (8.00e-03)} &\tabincell{c}{0.98 \\ (1.42e-05)} \\ \hline

  \tabincell{c}{Exam. 36 McCormick\\ Function \cite{SB2013} \\ (n = 12, m = n/2)}   & \tabincell{c}{12 \\ (3.00e-03)} &\tabincell{c}{1.31 \\ (1.05e-07)} & \tabincell{c}{13 \\ (1.00e-04)} &\tabincell{c}{1.31 \\ (8.42e-07)} & \tabincell{c}{5 \\ (5.00e-03)} &\tabincell{c}{1.31 \\ (2.64e-09)} \\ \hline

  \tabincell{c}{Exam. 37 Perm\\ Function $d, \ \beta$ \cite{SB2013} \\ (n = 4, m = n/2)}   & \tabincell{c}{25 \\ (4.00e-03)} &\tabincell{c}{1.19e+03 \\ (2.59e-06)} & \tabincell{c}{40 \\ (1.70e-02)} &\tabincell{c}{1.19e+03 \\ (3.85e-06)} & \tabincell{c}{28 \\ (8.00e-03)} &\tabincell{c}{1.19e+03 \\ (7.41e-05)} \\ \hline

  \tabincell{c}{Exam. 38 Power\\Sum Function \cite{SB2013} \\ (n = 4, m = n/2)}   & \tabincell{c}{1 \\ (6.00e-03)} &\tabincell{c}{1.52e+04 \\ (2.59e-06)} & \tabincell{c}{1 \\ (2.00e-03)} &\tabincell{c}{1.52e+04 \\ (1.85e-13)} & \tabincell{c}{2 \\ (6.00e-03)} &\tabincell{c}{1.52e+04 \\ (1.14e-12)} \\ \hline

  \tabincell{c}{Exam. 39 Price\\ Function \cite{AD2005} \\ (n = 2, m = n/2)}   & \tabincell{c}{8 \\ (6.00e-03)} &\tabincell{c}{7.06 \\ (1.50e-07)} & \tabincell{c}{9 \\ (1.10e-02)} &\tabincell{c}{7.06 \\ (1.85e-07)} & \tabincell{c}{11 \\ (6.00e-03)} &\tabincell{c}{7.06 \\ (7.77e-06)} \\ \hline

  \tabincell{c}{Exam. 40 Bohachevsky\\ Function \cite{SB2013} \\ (n = 2, m = n/2)}   & \tabincell{c}{9 \\ (4.00e-03)} &\tabincell{c}{2.36 \\ (1.75e-07)} & \tabincell{c}{11 \\ (4.00e-03)} &\tabincell{c}{2.36 \\ (7.12e-07)} & \tabincell{c}{9 \\ (6.00e-03)} &\tabincell{c}{1.15 \\ (1.00e-08)} \\ \hline

  \tabincell{c}{Exam. 41 Colville\\ Function \cite{SB2013} \\ (n = 4, m = n/2)}   & \tabincell{c}{13 \\ (2.00e-03)} &\tabincell{c}{21.14 \\ (1.35e-07)} & \tabincell{c}{26 \\ (1.70e-02)} &\tabincell{c}{21.14 \\ (5.81e-07)} & \tabincell{c}{11 \\ (5.00e-03)} &\tabincell{c}{5.19 \\ (3.37e-05)} \\ \hline

  \tabincell{c}{Exam. 42 Drop\\ Wave Function \cite{SB2013} \\ (n = 2, m = n/2)}   & \tabincell{c}{10 \\ (2.00e-03)} &\tabincell{c}{-0.79 \\ (2.12e-07)} & \tabincell{c}{9 \\ (2.00e-03)} &\tabincell{c}{-0.7858 \\ (9.15e-07)} & \tabincell{c}{5 \\ (5.00e-03)} &\tabincell{c}{-0.79 \\ (3.73e-09)} \\ \hline

  \tabincell{c}{Exam. 43 Schaffer\\ Function \cite{SB2013} \\ (n = 2, m = n/2)}   & \tabincell{c}{14 \\ (5.00e-03)} &\tabincell{c}{0.61 \\ (7.18e-08)} & \tabincell{c}{13 \\ (3.00e-03)} &\tabincell{c}{0.61 \\ (1.04e-07)} & \tabincell{c}{8 \\ (5.00e-03)} &\tabincell{c}{0.61 \\ (1.62e-07)} \\ \hline

  \tabincell{c}{Exam. 44 Six-Hump\\Camel Function \cite{SB2013} \\ (n = 2, m = n/2)}   & \tabincell{c}{10 \\ (4.00e-03)} &\tabincell{c}{0.74 \\ (1.73e-07)} & \tabincell{c}{18 \\ (2.00e-03)} &\tabincell{c}{0.74 \\ (4.09e-07)} & \tabincell{c}{11 \\ (5.00e-03)} &\tabincell{c}{0.74 \\ (2.00e-08)} \\ \hline

  \tabincell{c}{Exam. 45 Three-Hump\\Camel Function \cite{SB2013} \\ (n = 2, m = n/2)}   & \tabincell{c}{15 \\ (4.00e-03)} &\tabincell{c}{0.55 \\ (3.91e-07)} & \tabincell{c}{24 \\ (1.00e-04)} &\tabincell{c}{0.55 \\ (4.82e-07)} & \tabincell{c}{7 \\ (4.00e-03)} &\tabincell{c}{0.55 \\ (1.46e-07)} \\ \hline

  \tabincell{c}{Exam. 46 Trecanni\\ Function \cite{AD2005} \\ (n = 2, m = n/2)}   & \tabincell{c}{11 \\ (4.00e-03)} &\tabincell{c}{2.36 \\ (2.78e-08)} & \tabincell{c}{13 \\ (1.00e-04)} &\tabincell{c}{2.36 \\ (4.91e-07)} & \tabincell{c}{9 \\ (4.00e-03)} &\tabincell{c}{2.36 \\ (4.19e-08)} \\ \hline

  \tabincell{c}{Exam. 47 Box Bettes\\Exponential Quadratic \\ Function \cite{AD2005} \\ (n = 3, m = 2)}   & \tabincell{c}{20 \\ (2.60e-02)} &\tabincell{c}{1.42e-11 \\ (9.96e-07)} & \tabincell{c}{35 \\ (8.00e-03)} &\tabincell{c}{3.01e-13 \\ (7.07e-07)} & \tabincell{c}{13 \\ (6.00e-03)} &\tabincell{c}{3.26e-16 \\ (2.57e-08)} \\ \hline

  \tabincell{c}{Exam. 48 Chichinad \\ Function \cite{AD2005} \\ (n = 2, m = n/2)}   & \tabincell{c}{8 \\ (5.00e-03)} &\tabincell{c}{8.01 \\ (9.15e-09)} & \tabincell{c}{9 \\ (5.00e-03)} &\tabincell{c}{8.01 \\ (4.47e-08)} & \tabincell{c}{6 \\ (3.00e-03)} &\tabincell{c}{-20.06 \\ (8.52e-07)} \\ \hline

  \tabincell{c}{Exam. 49 Eggholder\\ Function \cite{SB2013} \\ (n = 2, m = n/2)}   & \tabincell{c}{17 \\ (3.00e-03)} &\tabincell{c}{-69.16 \\ (4.17e-07)} & \tabincell{c}{22 \\ (2.00e-02)} &\tabincell{c}{-69.60 \\ (9.45e-07)} & \tabincell{c}{9 \\ (6.00e-03)} &\tabincell{c}{-69.16 \\ (5.77e-08)} \\ \hline

  \tabincell{c}{Exam. 50 Exp2 \\ Function \cite{AD2005} \\ (n = 2, m = n/2)}   & \tabincell{c}{11 \\ (2.00e-03)} &\tabincell{c}{9.19 \\ (1.53e-08)} & \tabincell{c}{15 \\ (3.00e-03)} &\tabincell{c}{8.45 \\ (6.63e-07)} & \tabincell{c}{6 \\ (5.00e-03)} &\tabincell{c}{9.19 \\ (6.60e-07)} \\ \hline

  \tabincell{c}{Exam. 51 Hansen \\ Function \cite{AD2005} \\ (n = 2, m = n/2)}   & \tabincell{c}{9 \\ (6.00e-03)} &\tabincell{c}{-12.10 \\ (3.98e-08)} & \tabincell{c}{8 \\ (2.00e-03)} &\tabincell{c}{-12.10 \\ (1.69e-07)} & \tabincell{c}{6 \\ (5.00e-03)} &\tabincell{c}{-32.36 \\ (5.21e-05)} \\ \hline

  \tabincell{c}{Exam. 52 Hartmann \\ 3-D Function \cite{SB2013} \\ (n = 3, m = 2)}   & \tabincell{c}{13 \\ (2.00e-03)} &\tabincell{c}{-3.84 \\ (1.36e-07)} & \tabincell{c}{22 \\ (2.00e-03)} &\tabincell{c}{-3.84 \\ (8.13e-07)} & \tabincell{c}{2 \\ (4.00e-03)} &\tabincell{c}{-1.31e-30 \\ (2.22e-16)} \\ \hline

  \tabincell{c}{Exam. 53 Holder \\ Table Function \cite{SB2013} \\ (n = 2, m = n/2)}   & \tabincell{c}{13 \\ (2.00e-03)} &\tabincell{c}{-1.68 \\ (1.65e-07)} & \tabincell{c}{16 \\ (2.00e-03)} &\tabincell{c}{-1.68 \\ (5.82e-07)} & \tabincell{c}{5 \\ (5.00e-03)} &\tabincell{c}{-3.51e-02 \\ (9.62e-07)} \\ \hline

  \tabincell{c}{Exam. 54 Michalewicz \\ Function \cite{SB2013} \\ (n = 2, m = n/2)}   & \tabincell{c}{15 \\ (4.00e-03)} &\tabincell{c}{-1.00 \\ (7.24e-08)} & \tabincell{c}{16 \\ (3.00e-03)} &\tabincell{c}{-1.00 \\ (6.12e-07)} & \tabincell{c}{2 \\ (4.00e-03)} &\tabincell{c}{-7.18e-12 \\ (1.20e-10)} \\ \hline

  \tabincell{c}{Exam. 55 Schaffer \\ Function N. 4 \cite{SB2013} \\ (n = 4, m = n/2)}   & \tabincell{c}{8 \\ (5.00e-03)} &\tabincell{c}{0.30 \\ (9.11e-07)} & \tabincell{c}{11 \\ (1.00e-04)} &\tabincell{c}{0.30 \\ (8.35e-07)} & \tabincell{c}{6 \\ (6.00e-03)} &\tabincell{c}{0.29 \\ (7.68e-07)} \\ \hline

  \tabincell{c}{Exam. 56 Trefethen \\4 Function \cite{AD2005} \\ (n = 2, m = n/2)}   & \tabincell{c}{9 \\ (5.00e-03)} &\tabincell{c}{1.20 \\ (7.61e-08)} & \tabincell{c}{24 \\ (3.00e-04)} &\tabincell{c}{-1.36 \\ (1.51e-07)} & \tabincell{c}{8 \\ (6.00e-03)} &\tabincell{c}{-2.02 \\ (7.45e-05)} \\ \hline

  \tabincell{c}{Exam. 57 Zettl \\ Function \cite{AD2005} \\ (n = 2, m = n/2)}   & \tabincell{c}{11 \\ (4.00e-03)} &\tabincell{c}{0.14 \\ (3.26e-08)} & \tabincell{c}{17 \\ (1.00e-04)} &\tabincell{c}{0.14 \\ (8.27e-07)} & \tabincell{c}{11 \\ (4.00e-03)} &\tabincell{c}{0.14 \\ (1.84e-09)} \\ \hline

\end{tabular}
\end{table}

\vskip 2mm

\vskip 2mm

\vskip 2mm

\section{Conclusions}

\vskip 2mm

In this paper, we give the regularization continuation method with the 
trust-region updating strategy (Rcmtr) for linearly constrained optimization 
problems. Moreover, we reveals and utilizes the linear conservation law of the
regularization method and the quasi-Newton method such that it does not
need to compute the correction step other than the previous continuation method.
The new continuation method uses the inverse of the regularization two-sided
projection of the Lagrangian Hessian as the pre-conditioner to improve its 
robustness, which is other than the previous quasi-Newton methods. Numerical results show that Rcmtr
is more robust and faster than the traditional optimization method such as SQP
(the built-in subroutine fmincon.m of the MATLAB2020a environment \cite{MATLAB}),
the recent continuation method such as Ptctr \cite{LLS2021} and the alternating
direction method of multipliers (ADMM \cite{BPCE2011}). Therefore, Rcmtr is worth
exploring further, and we will extend it to the nonlinearly constrained optimization
problem in the future.

\section*{Acknowledgments} This work was supported in part by Grant 61876199 from National
Natural Science Foundation of China, Grant YBWL2011085 from Huawei Technologies
Co., Ltd., and Grant YJCB2011003HI from the Innovation Research Program of Huawei
Technologies Co., Ltd.. The authors are grateful to the anonymous referee for
his comments and suggestions which greatly improve presentation of this paper.

\vskip 2mm

\noindent \textbf{Conflicts of interest/Competing interests:} Not applicable.

\vskip 2mm

\noindent \textbf{Availability of data and material (data transparency):} If it is requested, we will
provide the test data.

\vskip 2mm

\noindent \textbf{Code availability (software application or custom code):} If it is requested, we will
provide the code.

\end{document}